\documentstyle{article}
\begin{document}

\begin{center}
{\large \bf Hilbert C*-modules and related subjects \\
-- a guided reference overview}  \\
Michael Frank \copyright \\
frank@mathematik.uni-leipzig.d400.de or \\
frank@rz.uni-leipzig.de
\end{center}

\noindent
\rule{6in}{0.1mm}

\noindent
last update: 1.4.96

\vspace{1cm}

\noindent
{\bf \S 1 About}

\medskip \noindent
Hilbert C*-modules are a often used tool in operator theory and in theory of
operator algebras. Beside this, the theory of Hilbert C*-modules is very 
interesting on it's own. Interacting with the theory of operator algebras and
including ideas of non-commutative geometry it progresses and produces results
and new problems attracting attention.

\noindent
At the contrary, the pieces of the theory of Hilbert C*-modules are rather
scattered through the literature. Most publications explain only as many 
definitions and results as necessary for the striven for applications of this
theory. However, there are some papers and chapters in monographs collecting
parts of the theory as well as E.~C.~Lance's well-written lecture notes.

\noindent
The purpose of the present reference overview is to show a practicable way for
systematic studies of the theory of Hilbert C*-modules.
Great emphasis is put on the historical consistency of the presented sources
following the line of ideas and applications. Since the term ''Hilbert ...
modules'' is in use for at least five mathematically more or less different
concepts one has always to pay attention what kind of theory is considered.
For the convenience of the reader we list the basic publications for all
known concepts wherein the notion ''Hilbert ... modules'' appears.
As a guide we refer to some basic publications on Hilbert C*-modules
representing essential achievements of the theory.
A second guide gives a short list of research fields wherein Hilbert
C*-modules are in use very actively, and some publications representing
these ways of application.

\noindent
The reader has to take into account that the choice of the sources is dominated
by the author's research interests and linguistic proficiency, as well as by
the local availability of sources. He apologizes for a probable insufficient
representation of the work of some colleagues in the present overview. All
suggestions, corrections and supplements are welcome.
 
\vspace{1cm}

\noindent
{\bf \S 2 Guide (part I)}

\medskip \noindent
{\sl Roots of the quite different notions of ''Hilbert ... modules'':}

\smallskip \noindent
{\sc I.~Kaplansky}, 1953, \cite{Kaplansky:53},
{\sc H.~Widom}, 1956, \cite{Widom:56} :
AW*-algebras, inner product AW*-modules (Kaplansky-Hilbert modules).

\smallskip \noindent
{\sc R.~M.~Loynes}, 1965, \cite{Loynes:65}:
VH-spaces, LVH-spaces.

\smallskip \noindent
{\sc R.~G.~Swan}, 1962, \cite{Swan:62},
{\sc J.~Dixmier} and {\sc A.~Douady}, 1963, \cite{Dixmier/Douady:63}:
vector bundles, projective modules. \newline
{\sc A.~O.~Takahashi}, 1971, \cite{Takahashi:79/2,Takahashi:79/1,Takahashi:71},
{\sc K.~H.~Hofmann}, 1972, \cite{Hofmann:72},
{\sc M.~J.~Dupr\'e}, 1972, \cite{Dupre:72,Dupre:74,Dupre:76},
({\sc H.~Takemoto}, 1973-76, \cite{Takemoto:73/1,Takemoto:75,Takemoto:76},)
{\sc J.~Varela}, 1974, \cite{Varela:74}:
Hilbert bundles, continuous fields of Hilbert spaces and of Banach algebras /
A categorial equivalence between (F)Hilbert bundles on compact spaces $K$ and
Hilbert C($K$)-modules.

\smallskip \noindent
{\sc N.~Wiener}, {\sc P.~R.~Masani}, 1957-66,
        \cite{Wiener/Masani:57,Wiener/Masani:58,Wiener/Masani:60,Masani:66},
{\sc H.~H.~Goldstine} and {\sc L.~P.~Horwitz}, 1966,
\cite{Goldstine/Horwitz:66},
{\sc P.~P.~Saworotnow}, 1968, \cite{Saworotnow:68}:
Hilbert (H*-)modules over matrix algebras/Hilbert*-algebras. (Hilbert
H*-modules are Hilbert C*-modules iff the H*-algebra is finite dimensional,
i.e.~a matrix algebra.) [see the second list of references]

\smallskip \noindent
{\sc D.~Bures}, 1971, \cite{Bures:71}:
special W*-valued inner products on von Neumann algebras.

\smallskip \noindent
{\sc W.~L.~Paschke}, 1972/73, \cite{Paschke:72,Paschke:73}:
one {\bf trailblazing paper} in Hilbert C*-module theory.

\smallskip \noindent
{\sc M.~A.~Rieffel}, 1974, \cite{Rieffel:74/2,Rieffel:74/1}:
the other {\bf trailblazing papers}, about Hilbert C*-modules and (strong)
Morita equivalence of C*-algebras.

\smallskip \noindent
{\sc Y.~Kakihara}, 1979-1984,
\cite{Kakihara:80,Kakihara:82,Kakihara:83/1,Kakihara:83/2,Kakihara:84}
\cite{Kakihara/Teresaki:79}:
Hilbert B(H)-modules with trace class valued inner product.

\smallskip \noindent
{\sc J.~Pincket}, 1986, \cite{Pincket:86}:
inner product C*-modules where the values of the inner product belong to the
duals of the underlying C*-algebras.

\smallskip \noindent
({\sc A.~Frydryszak, L.~Jakobczyk}, 1988, \cite{Frydryszak/Jakobczyk:88} :
Hilbert modules over infinite-dimensional Grass\-man-Banach algebras.)

\smallskip \noindent
({\sc R.~G.~Douglas, V.~I.~Paulsen}, 1989, \cite{Douglas/Paulsen:89}
(\cite{Muhly/Solel:95}):
''Hilbert modules'' := Hilbert spaces with a special C*-module structure on
them, but without {\it C*-valued} inner product) [see the third list of
references]

\smallskip \noindent
({\sc N.~C.~Phillips}, 1989, \cite{Phillips:89/1}: Hilbert modules over
pro-C*-algebras (i.e., over inverse limits of C*-algebras))

\smallskip \noindent
{\sc G.~Zeller-Meier}, 1991, \cite{Zeller--Meier:91/2,Zeller--Meier:91/1}:
Banach-C*-modules equipped with an C*-algebra valued inner product which is
not necessarily C*-linear/C*-antilinear in its arguments.

\smallskip \noindent
{\sc D.~P.~Blecher}, 1995, \cite{Blecher:95/1,Blecher:96/1}:
Hilbert modules over non-self-adjoint operator algebras, extending the
concept of Hilbert C*-modules / an alternative approach to Hilbert
C*-modules.

\bigskip 
{\sl Useful papers about Hilbert C*-modules from an axiomatic point of view on 
the theory}:

\smallskip \noindent
{\sc I.~Kaplansky} \cite{Kaplansky:53}/
{\sc H.~Widom} \cite{Widom:56}/ 
{\sc W.~L.~Paschke} \cite{Paschke:73,Paschke:76}/
{\sc M.~A.~Rieffel} \cite{Rieffel:74/2,Rieffel:74/1}/
{\sc G.~G.~Kasparov} \cite{Kasparov:80}/
{\sc M.~J.~Dupr\'e, P.~A.~Fillmore} \cite{Dupre/Fillmore:81}/
{\sc M.~I.~Gekhtman} \cite{Gekhtman:84}/
{\sc E.~V.~Tro\-{\H{\i}}tsky} \cite{Troitskij:86}/
{\sc J.~Cuntz, N.~Higson} \cite{Cuntz/Higson:87}/
{\sc O.~G.~Filippov} \cite{Filippov:87}/
{\sc J.-F.~Havet} \cite{Havet:88}/
{\sc H.~Lin} \cite{Lin:91/3,Lin:91/2,Lin:92}/
{\sc G.~Zeller-Meier} \cite{Zeller--Meier:91/2,Zeller--Meier:91/1}/
{\sc S.~Zhang}, \cite{Zhang:91/1}/
{\sc M.~Hamana} \cite{Hamana:92}/
{\sc E.~C.~Lance} \cite{Lance:95}/
{\sc L.~G.~Brown, J.~A.~Mingo and Nien-Tsu Shen}, \cite{Brown/Mingo/Shen:94}/
{\sc V.~M.~Manuilov} \cite{Manujlov:94,Manujlov:94/2,Manujlov:94/3}/
{\sc M.~Frank} \cite{Frank:90/1,Frank:93/5,Frank:95/4}.

\bigskip 
{\sl Ph.D. thesises using Hilbert C*-modules essentially}:

\smallskip \noindent
{\sc A.~O.~Takahashi}, 1971, \cite{Takahashi:71}/
{\sc M.~J.~Dupr\'e}, 1972, \cite{Dupre:72}/
{\sc W.~L.~Paschke}, 1972, \cite{Paschke:72}/
{\sc W.~Beer}, 1981, \cite{Beer:81}/
{\sc J.~A.~Mingo}, 1982, \cite{Mingo:82/0}/
{\sc N.-T.~Shen}, 1982, \cite{Shen:82}/
{\sc V.~A.~Trofimov}, 1987, \cite{Trofimov:87/4}/
{\sc J.~Weidner}, 1987, \cite{Weidner:87}/ 
{\sc Y.~Yang}, 1987, \cite{Yang:87}/
{\sc M.~Lesch}, 1988, \cite{Lesch:88}/ 
{\sc S.~Zhang}, 1988, \cite{Zhang:88}/
{\sc M.~Frank}, 1988, \cite{Frank:88}/
{\sc B.~Abadie}, 1992, \cite{Abadie:92}/
{\sc Huu~Hung Bui}, 1992, \cite{Bui:92}/
{\sc S. P. Kaliszewski}, 1994, \cite{Kaliszewski:94}.

\bigskip 
{\sl Books / Chapters in books and monographs about Hilbert C*-modules}:

\smallskip \noindent
Book: {\sc E.~C.~Lance}, 1993, \cite{Lance:95}

\smallskip \noindent
Chapters: {\sc A.~S.~Mishchenko}, 1984, \cite{Mishchenko:84/2} /
{\sc B.~Blackadar}, 1986, \cite{Blackadar:86} /
{\sc V.~I.~Istr\v{a}\c{t}escu}, 1987, \cite{Istratescu:87} /
{\sc N.~C.~Phillips}, 1989, \cite{Phillips:89} /
{\sc K.~K.~Jensen/K.~Thomsen}, 1991, \cite{Jensen/Thomsen:91}
(\cite{Thomsen:88}) /
{\sc N.~E.~Wegge--Olsen}, 1993, \cite{Wegge--Olsen:93}
(\cite{Wegge--Olsen:89}) /
{\sc H.~Schr\"oder}, 1993, \cite{Schroeder:93}/
{\sc A.~Connes}, 1994, \cite{Connes:94}.

\vspace{1cm}

\noindent
{\bf \S 3 Guide (part II)}

\medskip \noindent
The aim of this guide is to claim some bigger areas of research, where
Hilbert C*-modules have been used very actively and successfully. Listing
an author's publication means to list it as an example of a useful application
of the theory of Hilbert C*-modules in that area of research. Of course,
the main contributions to the subjects listed below often rest on quite
different mathematical methods. To keep the guide short and instructive
we mention only a few publications. For further sources the reader has to
consult these publications and the references therein.

\bigskip   
{\sl K-theory and KK-theory of operator algebras (G. G. Kasparov's approach):}

\smallskip \noindent
{\sc G.~G.~Kasparov}, 1980-90,
                     \cite{Kasparov:82,Kasparov:83,Kasparov:88,Kasparov:85}/
{\sc M.~V.~Pimsner, D.~Voiculescu}, 1980, \cite{Pimsner/Voiculescu:80}/
{\sc G.~Skandalis}, 1984/88/91, \cite{Skandalis:84,Skandalis:88,Skandalis:91}/
{\sc W.~L.~Paschke}, 1985, \cite{Paschke:85}/
{\sc M.~V.~Pimsner}, 1985, \cite{Pimsner:85}/
{\sc E.~V.~Tro\H{\i}tsky}, 1985, \cite{Troitskij:85/3}/
{\sc B.~Blackadar}, 1986, \cite{Blackadar:86}/
{\sc J.~Cuntz}, 1986, \cite{Cuntz:86}/
{\sc N.~C.~Phillips}, 1987/89, \cite{Phillips:87,Phillips:89}/
{\sc M.~A.~Rieffel}, 1987, \cite{Rieffel:87/1}/
{\sc J.~A.~Packer}, 1988, \cite{Packer:88}/
{\sc J.~Rosenberg}, 1990, \cite{Rosenberg:90}/
{\sc K.~K.~Jensen, K.~Thomsen}, 1991, \cite{Jensen/Thomsen:91}/
{\sc S.~Zhang}, 1991/92/93, \cite{Zhang:91/4,Zhang:91/3,Zhang:92/4,Zhang:93/3}/
{\sc N.~E.~Wegge-Olsen}, 1993, \cite{Wegge--Olsen:93}/
and others.

\bigskip 
{\sl Strong Morita equivalence of C*-algebras and its application to group 
representation theory:}

\smallskip \noindent
{\sc M.~A.~Rieffel}, 1974, \cite{Rieffel:74/2,Rieffel:74/1}/
{\sc L.~G.~Brown, P.~Green, M.~A.~Rieffel}, 1977,
                              \cite{Brown/Green/Rieffel:77}/
{\sc P.~Green}, 1978, \cite{Green:78}/
{\sc W.~Beer}, 1981/82, \cite{Beer:81,Beer:82}/
{\sc H.~H.~Zettl}, 1982, \cite{Zettl:82/1,Zettl:82/2}/
{\sc F.~Combes, H.~H.~Zettl}, 1983, \cite{Combes/Zettl:83}/
{\sc F.~Combes}, 1984, \cite{Combes:84}/
{\sc I.~Putnam}, 1985, \cite{Putnam:85}/
{\sc R.~J.~Plymen}, 1986/90, \cite{Plymen:86,Plymen:90}/
{\sc T.~Kajiwara}, 1987, \cite{Kajiwara:87}/
{\sc J.~A.~Packer}, 1988, \cite{Packer:88}/
{\sc P.~Xu}, 1991, \cite{Xu:91}/
{\sc K.~Mansfield}, 1991, \cite{Mansfield:91}/
{\sc S.~P.~Kaliszewski}, 1994, \cite{Kaliszewski:95,Kaliszewski:94}/
{\sc S.~Echterhoff}, 1993-95, \cite{Echterhoff:94,Echterhoff:95}/
{\sc J.~Quigg and I.~Raeburn}, 1995,
\cite{Kaliszewski/Quigg/Raeburn:95,Quigg:95,Echterhoff/Raeburn:96/2}/
and others.

\bigskip   
{\sl Normal operator-valued weights (resp., conditional expectations) of
finite index between C*-algebras / Correspondences of C*-algebras:}

\smallskip \noindent
{\sc D.~Bures}, 1971, \cite{Bures:71}/
{\sc A.~Connes}, 1980, \cite{Connes:80/2}/
{\sc M.~V.~Pimsner, S.~Popa}, 1986, \cite{Pimsner/Popa:86}/
{\sc M.~Baillet, Y.~Denizeau, J.-F.~Havet}, 1988,
                                          \cite{Baillet/Denizeau/Havet:88}/
{\sc Y.~Watatani}, 1990, \cite{Watatani:90}/
{\sc M.~Frank}, 1993, \cite{Frank:93/1}/
and others.

\bigskip    
{\sl AW*-algebras and monotone complete C*-algebras:}

\smallskip \noindent
{\sc I. Kaplansky}, 1953, \cite{Kaplansky:53}/
{\sc H. Widom}, 1956, \cite{Widom:56}/
{\sc J.~D.~M.~Wright}, 1969, \cite{Wright:69/2}/
{\sc C.~Sunouchi}, 1971, \cite{Sunouchi:71}/
{\sc K.~Sait{\^o}}, 1971-..., \cite{Saito:71,Saito:74,Saito:78,Saito:79}/
{\sc H.~Takemoto}, 1973, \cite{Takemoto:73/1}/
{\sc E.~Azoff}, 1978, \cite{Azoff:78}/
{\sc O.~Takenouchi}, 1978, \cite{Takenouchi:78}/
{\sc M.~Hamana}, 1979-...,
\cite{Hamana:79/2,Hamana:79/1,Hamana:81,Hamana:82/2,Hamana:82/3,Hamana:92}/
{\sc M.~Ozawa}, 1980-..., \cite{Ozawa:84,Ozawa:85/1,Ozawa:85/2,Ozawa:86}/
{\sc G.~A.~Elliott, K.~Sait{\^o}, J.~D.~M.~Wright}, 1983,
                                            \cite{Elliott/Saito/Wright:83}/
{\sc G.~K.~Pe\-dersen}, 1984/86, \cite{Pedersen:84,Pedersen:86}/
{\sc N.~Azarnia}, 1985, \cite{Azarnia:85}/
{\sc M.~Frank}, 1991/93, \cite{Frank:93/1,Frank:95/4}/
and others.

\bigskip  
{\sl Completely positive mappings between C*-algebras:}

\smallskip \noindent
{\sc W.~L.~Paschke}, 1973, \cite{Paschke:73}/
{\sc C.~Anantharaman-Delaroche}, 1990, \cite{Anantharaman:90/2}/
{\sc C.~Anantharaman-Delaroche, J.-F.~Havet}, 1990,
                     \cite{Anantharaman/Havet:90}/
{\sc J.~A.~Mingo}, 1990, \cite{Mingo:90}/
{\sc H.~Lin}, 1991, \cite{Lin:91/3}/
{\sc J.~Tsui}, 1996, \cite{Tsui:96}/
and others.

\bigskip  
{\sl Mathematical and theoretical physics:}

\smallskip \noindent
{\sc M.~Banai}, 1987, \cite{Banai:87} /
{\sc D.~Applebaum}, 1988, \cite{Applebaum:87}/
{\sc P.~Xu}, 1991/92, \cite{Xu:91,Xu:92}/
{\sc N.~P.~Landsman}, 1993-94, \cite{Landsman:94,Landsman:95}/
{\sc V.~M.~Manuilov}, 1994, \cite{Manujlov:94/2}/
{\sc Yun-Gang Lu}, 1995, \cite{Lu:95}/
and others.

\bigskip  
{\sl Unbounded operators and quantum groups:}

\smallskip \noindent
{\sc M.~Hilsum}, 1989, \cite{Hilsum:89}/
{\sc S.~Baaj, P.~Julg}, 1983, \cite{Baaj/Julg:83}/
{\sc S.~Baaj, G.~Skandalis}, 1989, \cite{Baaj/Skandalis:89}/
{\sc S.~L.~Woronowicz}, 1991, \cite{Woronowicz:91}
{\sc S.~L.~Woronowicz, K.~Napi\'orkowski}, 1992,
                                       \cite{Woronowicz/Napiorkowski:92}/
{\sc E.~C.~Lance}, 1994, \cite{Lance:94}/
and others.

\bigskip   
{\sl Vector bundles, (F)Hilbert bundles $\leftrightarrow$ projective
C*-modules, Hilbert C*-modules:}

\smallskip \noindent
{\sc R.~G.~Swan}, 1962, \cite{Swan:62}/
{\sc J.~Dixmier, A.~Douady}, 1963, \cite{Dixmier/Douady:63}
{\sc A.~O.~Takahashi}, 1971, \cite{Takahashi:79/2,Takahashi:79/1,Takahashi:71}/
{\sc K.~H.~Hofmann}, 1972, \cite{Hofmann:72}/
{\sc M.~J.~Dupr\'e}, 1972, \cite{Dupre:72,Dupre:74,Dupre:76}/
{\sc J.~Varela}, 1974, \cite{Varela:74},
{\sc J.~Fell}, 1977, \cite{Fell:77}/
{\sc M.~A.~Rieffel}, 1983/85/88, \cite{Rieffel:83/2,Rieffel:85,Rieffel:88/1}/
{\sc A.~J.~L.~Sheu}, 1987, \cite{Sheu:87/2}/
{\sc R.~G.~Swan}, 1987, \cite{Swan:87}/
{\sc J.~M.~S.~Fell, R.~S.~Doran}, 1988, \cite{Fell/Doran:88}/
and others.

\bigskip  
{\sl Rotation C*-algebras and related structures:}

\smallskip \noindent
{\sc M.~A.~Rieffel}, 1981-83, \cite{Rieffel:81,Rieffel:83/2}/
{\sc M.~{De Brabander}}, 1984, \cite{DeBrabander:84}/
{\sc J.~A.~Packer}, 1987/88, \cite{Packer:87,Packer:88}/
{\sc S.~G.~Walters}, 1994/95, \cite{Walters:94,Walters:95}/
{\sc G.~A.~Elliott, Lin Qing}, 1995, \cite{Elliott/Qing:95,Qing:95}/
and others.

\newpage

\bigskip  
{\sl Other non-commutative geometry:}

\smallskip \noindent
{\sc A.~Connes}, 1980-...,
                 \cite{Connes:80,Connes:82,Connes:86,Connes:90/2,Connes:90}/
{\sc M.~A.~Rieffel}, 1990, \cite{Rieffel:90/1}/
{\sc J.~Lott, A.~Connes}, 1992, \cite{Lott/Connes:92}/
{\sc J.~Cuntz}, 1993, \cite{Cuntz:93}/
and others.

\bigskip   
{\sl Topological applications:}

\smallskip \noindent
{\sc G.~G.~Kasparov}, 1975-..., \cite{Kasparov:75,Kasparov:83,Kasparov:88}/
{\sc A.~S.~Mishchenko}, 1978-79, \cite{Mishchenko:78,Mishchenko:79/1}/
{\sc A.~S.~Mi\-shchenko, A.~T.~Fomenko}, 1979, \cite{Mishchenko/Fomenko:79/2}/
{\sc R.~A.~Biktashev, A.~S.~Mishchenko}, 1980, \cite{Biktashev/Mishchenko:80}/
{\sc V.~A.~Kasimov}, 1982, \cite{Kasimov:82/1}/
{\sc M.~A.~Rieffel}, 1983, \cite{Rieffel:83/1}/
{\sc J.~Rosenberg}, 1983-90, \cite{Rosenberg:83,Rosenberg:88,Rosenberg:90}/
{\sc M.~Hilsum, G.~Skandalis}, 1983/1992,
                      \cite{Hilsum/Skandalis:83,Hilsum/Skandalis:92}/
{\sc J.~Kaminker, J.~G.~Miller}, 1985, \cite{Kaminker/Miller:85}/
{\sc F.~Sharipov}, 1985, \cite{Sharipov:85/1}/
{\sc E.~V.~Tro\H{\i}tsky}, 1985,
\cite{Troitskij:85/4,Troitskij:85/1,Troitskij:85/2,Troitskij:88,Troitskij:89}/
{\sc G.~G.~Kasparov, G.~Skandalis}, 1990/91,
                     \cite{Kasparov/Skandalis:90,Kasparov/Skandalis:91}/
and others.

\bigskip    
{\sl Hilbert modules over pro-C*-algebras:}

\smallskip \noindent
{\sc N.~C.~Phillips}, 1989, \cite{Phillips:89/1}/
{\sc J.~Weidner}, 1989, \cite{Weidner:89/1,Weidner:89/2}/
{\sc C.~Schochet}, 1994, \cite{Schochet:94}/
and others.

\bigskip    
{\sl Prediction theory of multivariate stochastic processes:}

\smallskip \noindent
{\sc N.~Wiener}, {\sc P.~R.~Masani}, 1957-66,
               \cite{Wiener/Masani:57,Wiener/Masani:58,Wiener/Masani:60},
                 \cite{Masani:59,Masani:62,Masani:66}/
{\sc R.~M.~Loynes}, 1965, \cite{Loynes:65/3}/
{\sc H. Salehi}, 1965-67, \cite{Salehi:65,Salehi:66,Salehi:67}/
{\sc P.~P.~Saworotnow}, 1983, \cite{Saworotnow:83}/
{\sc H.~Fuge}, 1995, \cite{Fuge:95} /
{\sc U.~Gerecke}, {\sc J.~Lorenz}, 1995, \cite{Gerecke/Lorenz:95}/
{\sc A.~Kokschal}, 1995, \cite{Kokschal:95}/
and others.

\bigskip \noindent
{\bf Acknowledgement:} I wish to thank all the colleagues who submitted
preprints, reprints and copies of their Ph.D.'s, and who suggested sources
concerned with the subject.
Especially, I am indebted to E.~V.~Troitsky who brought his literature
list to my attention, and to B.~Kirstein who showed me the use of Hilbert
modules in stochastics.

\newpage

{\small

\nocite{Filippov:87,Filippov:90/1,Filippov:90/2,Kasimov:82/1}
\nocite{Irmatov:90,Irmatov:89,Mishchenko:79/1,Mishchenko:84}
\nocite{Mishchenko/Fomenko:79/2,Mishchenko/Sharipov:83}
\nocite{Trofimov:86,Trofimov:87/1,Trofimov:87/2}
\nocite{Trofimov:87/3,Trofimov:87/4,Troitskij:85/1,Troitskij:85/2}
\nocite{Troitskij:85/3,Anantharaman:87/2,Baaj/Julg:83,Brown/Pedersen:93}
\nocite{Connes:86,Frank:95/4,Frank:88,Hilsum/Skandalis:83}
\nocite{Dupre/Fillmore:81,Kasparov:75,Kasparov:80,Kasparov:81}
\nocite{Mingo:82,Paschke:73,Lin:91/3,Lin:92,Lin:91}
\nocite{Baillet/Denizeau/Havet:88,Mingo/Phillips:84,Brown/Green/Rieffel:77}
\nocite{Rieffel:74/2,Rieffel:83/2,Fox/Haskell/Raeburn:89}
\nocite{Anantharaman:90,Anantharaman/Havet:90,Anantharaman:87/1}
\nocite{Arkhangel'skij:88,Atiyah/Bott/Shapiro:64,Azarnia:85,Azoff:78}
\nocite{Baaj/Skandalis:89,Beer:82,Beer:81,Biktashev:83}
\nocite{Biktashev/Mishchenko:80,Blackadar:86,Blackadar:83,Brown:88}
\nocite{Combes/Zettl:83,Combes:84,Connes:80,Connes:82,Connes/Skandalis:84}
\nocite{Connes/Rieffel:87,Cuntz/Higson:87,Curto/Muhly/Williams:84,Dadarlat:88}
\nocite{Dixmier/Douady:63,Dupre/Gillette:83}
\nocite{Elliott/Saito/Wright:83,Elliott:84,Fell:77,Fell/Doran:88}
\nocite{Frank:85,Frank:90/1,Frank:89/1,Frank:89/2,Frank:90/2,Frank:90/3}
\nocite{Giordano/Handelman:89,Goldstine/Horwitz:66,Green:78,Greene:74}
\nocite{Haagerup:89,Halpern:85,Hennings:89,Havet:88,Higson:88}
\nocite{Hilsum:85,Hofmann:72,Hofmann:74,Itoh:80,Itoh:90}
\nocite{Kakihara:82,Kasimov:81,Kurmakaev/Shkarin:90,Paschke:85}
\nocite{Kajiwara:87,Kaminker/Miller:85,Kakihara/Teresaki:79,Kakihara:80}
\nocite{Kakihara:83/1,Kakihara:83/2,Kakihara:84,Kakihara:85,Kaplansky:53}
\nocite{Kasimov:82/2,Kasparov:83,Kasparov:85,Moore/Schochet:88}
\nocite{Kasparov:88,Kasparov/Skandalis:90,Kodaka:89/2,Kodaka:89/1}
\nocite{Kumijan:88,Lee/Kim:82,Loynes:65,Mallios:85,Mingo:87,Ozawa:84}
\nocite{Mishchenko:84/2,Mishchenko:79/3,Mishchenko/Filippov:89}
\nocite{Mori:88,Muhly/Renault/Williams:87,Nagisa/Song:89,Ozawa:80,Ozawa:83}
\nocite{Ozawa:90,Ozawa:85/1,Ozawa:85/2,Ozawa:86,Ozawa/Saito:86}
\nocite{Pimsner/Voiculescu:80,Rieffel:82/2}
\nocite{Packer:86,Packer:87,Packer:88,Paschke:74,Paschke:76,Paschke:77}
\nocite{Phillips:89,Phillips:88/1,Phillips:88/2,Phillips:87}
\nocite{Pimsner:85,Pincket:86,Plymen:86,Plymen:87,Putnam:85,Raeburn:88}
\nocite{Putnam/Schmidt/Skau:86,Rieffel:74/1,Rieffel:76,Rieffel:82/1}
\nocite{Rieffel:85,Rieffel:83/1,Rieffel:83/2,Rieffel:87/1,Rieffel:87/2}
\nocite{Rieffel:88/1,Rieffel:90/1,Rieffel:90/2,Rosenberg:88,Saito:71,Saito:74}
\nocite{Sauvageot:89,Saworotnow:68,Scedrow/Scowcroft:88,Sharipov:85/1}
\nocite{Sharipov:85/2,Sharipov/Zhuraev:88,Sheu:87/1,Sheu:87/2,Umegaki:55}
\nocite{Skandalis:84,Skandalis:88,Sunouchi:71,Swan:62,Takeuti:83/1}
\nocite{Swan:87,Swan:77,Takahashi:71,Takemoto:75,Takemoto:73/1,Takesaki:60}
\nocite{Takeuti:83/2,Thomsen:88,Tomiyama:87,Troitskij:85/4,Troitskij:86}
\nocite{Vincent--Smith:67,Wang:89,Watatani:90,Weidner:87,Widom:56,Wittstock:84}
\nocite{Wittstock:81,Wright:69/2,Yamagami:84,Zekri:89,Zeller--Meier:87}
\nocite{Zeller--Meier:91/2,Zettl:82/1,Zettl:82/2,Zettl:83,Mingo:90}
\nocite{Douglas/Paulsen:89,Wegge--Olsen:89,Zhang:93,Zeller--Meier:91/1}
\nocite{Shen:82,Woronowicz:91,Woronowicz/Napiorkowski:92,Yang:87,Nishimura:90}
\nocite{Applebaum:87,Takemoto:76,Rieffel:91,Connes:90,Takahashi:79/2}
\nocite{Higson:90,Rosenberg:90,Plymen:90,Lin:91/2,Takahashi:79/1}
\nocite{Kasparov/Skandalis:91,Frank:92/6,Connes:90/2,,Mishchenko:78}
\nocite{Bunke:92,Connes:90,Havet:90,Jensen/Thomsen:91,Mishchenko/Solov'ov:77}
\nocite{Arveson:91,Brown:85,Kaftal:91,Xu:91,Lott/Connes:92,Olsen/Pedersen:89}
\nocite{Packer/Raeburn:90,Packer/Raeburn:89,Saito:78,Saito:79,Skandalis:91}
\nocite{Takenouchi:78,Troitskij:87,Zhang:90,Zhang:91/1,Zhang:91/2,Zhang:92}
\nocite{Wegge--Olsen:93,Connes:81,Hilsum:89,Lesch:88,Exel:93,Zhang:91/4}
\nocite{Istratescu:87,Aubert:76,Bures:71,Zhang:92/2,Brown/Mingo/Shen:94}
\nocite{Abadie:92,Abadie:93/1,Abadie:94,Zhang:93/2,Dupre:72,Dupre:74,Dupre:76}
\nocite{Paschke:72,Hamana:92,Anantharaman:90/2,Tsui:96,Cuntz:93,Frank:93/1}
\nocite{Weidner:89/1,Weidner:89/2,Lance:94,Kasimov:89/1,Hamana:92}
\nocite{Kasimov:89/2,Loginov/Shulman:93,Denizeau/Havet:94,Denizeau/Havet:93/2}
\nocite{Jolissaint:91,Troitskij:86/3,Troitskij:88,Troitskij:89,Cuntz:86}
\nocite{Mishchenko:84/2,Packer:88/2,Sharipov/Zhuraev:86,Raeburn:81}
\nocite{Troitskij:86/2,Brown/Pedersen:93,Zhang:91/3,Zhang:88,Zhang:89}
\nocite{Zhang:91/5,Zhang:91/6,Zhang:93/3,Zhang:92/4,Biktashev:82,Je/Yang:86}
\nocite{Ghez/Lima/Roberts:85,Arzumanian/Grigorian:90,Troitskij:86/2}
\nocite{Kandelaki:86,Frank:93/5,Hegerfeldt:85,Lance:95,Raeburn/Williams:93}
\nocite{Landsman:95,Troitskij:93/2,Pimsner:93,Troitskij:94,DeBrabander:84}
\nocite{Parker:88,Giordano:88,Gekhtman:84,Carey/Phillips:91,Yang:84}
\nocite{Landsman:94,Rieffel:93,Papatriantafillou:94,Ji:94,Manujlov:94}
\nocite{Manujlov:94/2,Irmatov/Mishchenko:90,Bui:94,Bui:95/1,Bui:95/2}
\nocite{Blecher/Muhly/Paulsen:94,Lin:94,Putnam:88,Delanghe:76,Lin:93}
\nocite{Raeburn/Williams:85,Pimsner/Popa:86,Bultheel:82,Rieffel:81}
\nocite{Schochet:94,Rieffel:93/2,Kasparov:82,Manujlov:94/3,Troitskij:94/2}
\nocite{Troitskij:91,Yamagami:94,Kaliszewski:94,Exel:94,Lin:93/2,Xu:92}
\nocite{Manujlov:95/1,Abadie/Eilers/Exel:95,Frank/Troitsky:95,Kim/Yang:84}
\nocite{Manujlov:95/2,Echterhoff:90,Echterhoff:93,Echterhoff:94,Echterhoff:95}
\nocite{Echterhoff:94/2,Echterhoff/Raeburn:96/2,Troitskij:95,Bui:92}
\nocite{Jeong:93,Saito:94,Schroeder:93,Zhang:94/1,Weaver:95/1,Weaver:95/2}
\nocite{Frank:95/3,Quigg/Spielberg,Lutz:95,Troitskij:92,Kajiwara/Watatani:95}
\nocite{Phillips/Raeburn:94,Magajna:95/1,Magajna:95/2,Bui:94/4}
\nocite{Elliott/Qing:95,Qing:95,Blecher:95/1,Blecher:96/1,Magajna:94,Ng:95}
\nocite{Walters:94,Walters:95,Ng:95/2,Ng:95/3,Manujlov:94/4,Raeburn:95}
\nocite{Nishimura:93,Frydryszak/Jakobczyk:88,Bunke:95,Banai:87,Wickstead:82}
\nocite{Blecher:95/3,Muhly/Solel:95,Heo:96,Kirchberg:91,Mansfield:91}
\nocite{Quigg:95,Kaliszewski/Quigg:95,Kaliszewski/Quigg/Raeburn:95}
\nocite{Echterhoff/Raeburn:96/2,Kaliszewski:95,Blecher/Muhly/Na:96/1}
\nocite{Bui:95/3,Bui:96/1,Damir:94,Ng:96,Frank:96/3,Wiener/Masani:57}
\nocite{Wiener/Masani:58,Wiener/Masani:60,Masani:62,Masani:66,Masani:59}
\nocite{Kokschal:95,Fuge:95,Saworotnow:83,Loynes:65/2,Loynes:65/3}
\nocite{Truong-Van:81,Salehi:65,Salehi:67,Salehi:66,Troitskij:96}
\nocite{Rosenberg:64,Gerecke/Lorenz:95,Lu:95,Murphy:96,Manujlov:96}
\nocite{Khimshiashvili:92,Ellis/Gohberg/Lay:95,Ben-Artzi/Gohberg:94}
\nocite{Irmatov:88,Echterhoff:94/3,Baaj/Skandalis:93,Cuntz/Skandalis:86}
\nocite{Kasparov:93,Mingo:82/0,Blanchard:95,Hilsum/Skandalis:92}
\nocite{Zhang:92/3,Frank/Manuilov:95}



\newpage

\noindent
The graphic below lists the number of publications published in a year. 
The number for 1994, 1995 and 1996 contains some circulating preprints
the final year of publication of which has to be defined through the years.

\bigskip
\begin{tabular}{ccl}
year & number of publications & $mm$ \\[0.5ex]
1953 & 1 & $\rule{1mm}{3mm}$ \\
1954 & 0 & \\
1955 & 1 & $\rule{1mm}{3mm}$ \\
1956 & 1 & $\rule{1mm}{3mm}$ \\
1957 & 1 & $\rule{1mm}{3mm}$ \\
1958 & 1 & $\rule{1mm}{3mm}$ \\
1959 & 1 & $\rule{1mm}{3mm}$ \\
1960 & 2 & $\rule{2mm}{3mm}$ \\
1961 & 0 & \\
1962 & 2 & $\rule{2mm}{3mm}$ \\
1963 & 1 & $\rule{1mm}{3mm}$ \\
1964 & 2 & $\rule{2mm}{3mm}$ \\
1965 & 4 & $\rule{4mm}{3mm}$ \\
1966 & 3 & $\rule{3mm}{3mm}$ \\
1967 & 2 & $\rule{2mm}{3mm}$ \\
1968 & 1 & $\rule{1mm}{3mm}$ \\
1969 & 1 & $\rule{1mm}{3mm}$ \\
1970 & 0 & \\
1971 & 4 & $\rule{4mm}{3mm}$ \\
1972 & 3 & $\rule{3mm}{3mm}$ \\
1973 & 2 & $\rule{2mm}{3mm}$ \\
1974 & 8 & $\rule{8mm}{3mm}$ \\
1975 & 2 & $\rule{2mm}{3mm}$ \\
1976 & 6 & $\rule{6mm}{3mm}$ \\
1977 & 5 & $\rule{5mm}{3mm}$ \\
1978 & 5 & $\rule{5mm}{3mm}$ \\
1979 & 8 & $\rule{8mm}{3mm}$ \\
1980 & 11 & $\rule{11mm}{3mm}$ \\
1981 & 10 & $\rule{10mm}{3mm}$ \\
1982 & 15 & $\rule{15mm}{3mm}$ \\
1983 & 18 & $\rule{18mm}{3mm}$ \\
1984 & 18 & $\rule{18mm}{3mm}$ \\
1985 & 24 & $\rule{24mm}{3mm}$ \\
1986 & 19 & $\rule{19mm}{3mm}$ \\
1987 & 29 & $\rule{29mm}{3mm}$ \\
1988 & 32 & $\rule{32mm}{3mm}$ \\
1989 & 27 & $\rule{27mm}{3mm}$ \\
1990 & 29 & $\rule{29mm}{3mm}$ \\
1991 & 24 & $\rule{24mm}{3mm}$ \\
1992 & 18 & $\rule{18mm}{3mm}$ \\
1993 & 26 & $\rule{26mm}{3mm}$ \\
(1994) & 27 & $\rule{27mm}{3mm}$ \\
(1995) & 40 & $\rule{40mm}{3mm}$ \\
(1996) & 13 & $\rule{13mm}{3mm}$
\end{tabular}

}


\begin{thebibliography}{100}

\bibitem{Abadie:92}
B.~Abadie.
\newblock {\em {On the K-theory of non-commutative Heisenberg manifolds}}.
\newblock PhD thesis, Univ. of California at Berkeley, 1992.

\bibitem{Abadie:93/1}
B.~Abadie.
\newblock {Generalized fixed-point algebras of certain actions on crossed
  products}.
\newblock preprint, funct-an@babbage.sissa.it , \# 9301004, 1993.

\bibitem{Abadie:94}
B.~Abadie.
\newblock {''Vector bundles'' over quantum Heisenberg manifolds}.
\newblock In R.~Curto and P.~E.~T. J{\o}rgensen, editors, {\em Algebraic
  Methods in Operator Theory}. Birkh\"auser, Boston - Basel - Berlin, 1994.

\bibitem{Abadie/Eilers/Exel:95}
B.~Abadie, S.~Eilers, and R.~Exel.
\newblock {Morita equivalence for crossed products by Hilbert C*-bimodules}.
\newblock preprint, K{\o}benhavns Universitet, Matematisk Institut, Denmark,
  1995.

\bibitem{Anantharaman:87/1}
C.~Anantharaman-Delaroche.
\newblock {On Connes property T for von Neumann algebras}.
\newblock {\em Math. Japo\-nica}, {\bf 32}, (1987).
\newblock 337-355.

\bibitem{Anantharaman:87/2}
C.~Anantharaman-Delaroche.
\newblock {Syst\`emes dynamiques non commutatifs et moyennabilit\'e}.
\newblock {\em Math. Ann.}, {\bf 279}, (1987).
\newblock 297-315.

\bibitem{Anantharaman:90/2}
C.~Anantharaman-Delaroche.
\newblock {On completely positive maps defined by an irreducible
  correspondence}.
\newblock {\em Can. Math. Bull.}, {\bf 33}, (1990).
\newblock 434-441.

\bibitem{Anantharaman:90}
C.~Anantharaman-Delaroche.
\newblock {On relative amenability for von Neumann algebras}.
\newblock {\em Compos. Math.}, {\bf 74 }, (1990).
\newblock 333-352.

\bibitem{Anantharaman/Havet:90}
C.~Anantharaman-Delaroche and J.~F. Havet.
\newblock {On approximate factorizations of completely positive maps}.
\newblock {\em J. Funct. Anal.}, {\bf 90}, (1990).
\newblock 411-428.

\bibitem{Applebaum:87}
D.~Applebaum.
\newblock {Quantum stochastic parallel transport on non-commutative vector
  bundles}.
\newblock In {\em Quantum probability and applications, III, Oberwolfach 1987}.
  Lecture Note Math. {\bf 1303} , Springer--Verlag, Berlin, 1988.
\newblock pp. 20-36.

\bibitem{Arkhangel'skij:88}
A.~V. Arkhangel'skij.
\newblock {\em {General topology - 2}}.
\newblock {Itogi nauki i tekhniki, Sovrem. probl. mat., fund. naprawl., v. {\bf
  50}}. VINITI, Moscou, 1988.

\bibitem{Arveson:91}
W.~Arveson.
\newblock {C*-algebras associated with sets of semigroups of isometries}.
\newblock {\em Internat. J. Math.}, {\bf 2}, (1991).
\newblock 235-255.

\bibitem{Arzumanian/Grigorian:90}
V.~A. Arzumanian and S.~A. Grigorian.
\newblock {Invariant algebras of operator fields on compact abelian groups
  (russ./engl.)}.
\newblock {\em Izv. Akad. Nauk Armyan. SSR, Ser. Mat.}, {\bf 25}, (1990).
\newblock no. 4, 333-343 / {\it Soviet J. Contemp. Math. Anal.} {\bf 25}(1990),
  no. 4, 20-31.

\bibitem{Atiyah/Bott/Shapiro:64}
M.~F. Atiyah, R.~Bott, and A.~Shapiro.
\newblock {Clifford modules}.
\newblock {\em Topology}, {\bf 3}, (1964).
\newblock suppl. 1, 3-38.

\bibitem{Aubert:76}
P.-L. Aubert.
\newblock {Th\'eorie de Galois pour une W*-Alg\`ebre}.
\newblock {\em Comment. Math. Helvetici}, {\bf 51}, (1976).
\newblock 411-433.

\bibitem{Azarnia:85}
N.~Azarnia.
\newblock {Dense operator on a KH-module}.
\newblock {\em Rend. Circ. Mat. Palermo}, {\bf 34}, (1985).
\newblock 105-110.

\bibitem{Azoff:78}
E.~Azoff.
\newblock {Kaplansky--Hilbert modules and self-adjointness of operator
  algebras}.
\newblock {\em Amer. J. Math.}, {\bf 100}, (1978).
\newblock 957-972.

\bibitem{Baaj/Julg:83}
S.~Baaj and P.~Julg.
\newblock {Th\'eorie bivariante de Kasparov et op\'erateurs non born\'es dans
  les C*-modules hilbertiens}.
\newblock {\em C. R. Acad. Sci. Paris, s\'er. 1}, {\bf 296}, (1983).
\newblock 875-878.

\bibitem{Baaj/Skandalis:89}
S.~Baaj and G.~Skandalis.
\newblock {C*-alg\`ebres de Hopf et th\'eorie de Kasparov \'equivariante}.
\newblock {\em K-theory}, {\bf 2}, (1989).
\newblock 683-721.

\bibitem{Baaj/Skandalis:93}
S.~Baaj and G.~Skandalis.
\newblock {Unitaires multiplicatifs et dualit\'e pour les produits crois\'es de
  C*-alg\`ebres}.
\newblock {\em Ann. Sci. \'Ec. Norm. Sup., 4e s\'er.}, {\bf 26}, (1993).
\newblock 425-488.

\bibitem{Baillet/Denizeau/Havet:88}
M.~Baillet, Y.~Denizeau, and J.-F. Havet.
\newblock {Indice d'une esperance conditionelle}.
\newblock {\em Comp. Math.}, {\bf 66}, (1988).
\newblock 199-236.

\bibitem{Banai:87}
M.~Banai.
\newblock {An unconventional canonical quantization of local scalar fields
  over quantum space-time}.
\newblock {\em J. Math. Phys.}, {\bf 28}, (1987).
\newblock 193-214.

\bibitem{Beer:81}
W.~Beer.
\newblock {\em {On Morita equivalence of C*-algebras}}.
\newblock PhD thesis, Univ. of California, Berkeley, USA, 1981.

\bibitem{Beer:82}
W.~Beer.
\newblock {On Morita equivalence of nuclear C*-algebras}.
\newblock {\em J. Pure and Appl. Algebra}, {\bf 26}, (1982).
\newblock 249-267.

\bibitem{Ben-Artzi/Gohberg:94}
A.~Ben-Artzi and I.~Gohberg.
\newblock {Orthogonal polynomials over Hilbert modules}.
\newblock In A.~Feintuch et~al., editor, {\em {Nonselfadjoint operators and
  related topics. Workshop on Operator theory and its applications, Beersheva,
  Israel, February 24-28, 1992}}, volume~{\bf 73} of {\em Basel,
  Birkh\"auser-Verlag, {\it Oper. Theory, Adv. Appl.}}, (1994).
\newblock 96-126.

\bibitem{Biktashev:82}
R.~A. Biktashev.
\newblock {Spectra of pseudodifferential operators over C*-algebras
  (russ./engl.)}.
\newblock {\em Vestnik Moskov. Univ., Ser. I: Mat.-Mekh.}, no. 4, (1982).
\newblock 36-38 / {\it Moscow Univ. Math. Bull.} {\bf 37}(1982), no. 4, 45-48.

\bibitem{Biktashev:83}
R.~A. Biktashev.
\newblock {Spectra of pseudodifferential operators over C*-algebras (russ.)}.
\newblock {\em Trudy Semin. Vektor Tenzor Anal.}, {\bf 21}, (1983).
\newblock 259-267.

\bibitem{Biktashev/Mishchenko:80}
R.~A. Biktashev and A.~S. Mishchenko.
\newblock {Spectra of elliptic unbounded pseudodifferential operators over
  C*-algebras (russ./engl.)}.
\newblock {\em Vestn. Mosk. Univ., Ser. I: Mat.-Mekh.}, no. 3, (1980).
\newblock 56-58 / {\it Moscow Univ. Math. Bull.} {\bf 35}(1980), no. 3, 59-65.

\bibitem{Blackadar:83}
B.~Blackadar.
\newblock {A stable cancellation theorem for simple C*-algebras}.
\newblock {\em Proc. London Math. Soc.}, {\bf 47}, (1983).
\newblock 303-305.

\bibitem{Blackadar:86}
B.~Blackadar.
\newblock {\em {K-theory for operator algebras}}.
\newblock Springer--Verlag, New York, 1986.

\bibitem{Blanchard:95}
E.~Blanchard.
\newblock {Tensor products of C(X)-algebras over C(X)}.
\newblock {\em Ast\'erisque}, {\bf 232}, (1995).
\newblock 81-92.

\bibitem{Blecher:95/1}
D.~P. Blecher.
\newblock {A new approach to Hilbert C*-modules}.
\newblock preprint, University of Houston, Houston, Texas, U.S.A. / to appear
  in {\it Math. Ann.}, 1995.

\bibitem{Blecher:95/3}
D.~P. Blecher.
\newblock {On selfdual Hilbert modules}.
\newblock preprint, University of Houston, Houston, Texas, U.S.A., submitted to
  J. Funct. Anal., 1995.

\bibitem{Blecher:96/1}
D.~P. Blecher.
\newblock {A generalization of Hilbert modules}.
\newblock {\em J. Funct. Anal.}, {\bf 136}, (1996).
\newblock 365-421.

\bibitem{Blecher/Muhly/Na:96/1}
D.~P. Blecher, P.~S. Muhly, and Qiyuan Na.
\newblock {Morita equivalence of operator algebras and their C*-envelopes}.
\newblock preprint, 1996.

\bibitem{Blecher/Muhly/Paulsen:94}
D.~P. Blecher, P.~S. Muhly, and V.~I. Paulsen.
\newblock {Categories of operator modules -- Morita equivalence and projective
  modules}.
\newblock preprint, University of Iowa, Iowa City, U.S.A. / to appear in {\it
  Memoirs Amer. Math. Soc.}, 1994.

\bibitem{Brown/Pedersen:93}
L.~G. Brown.
\newblock {On the geometry of the unit ball of a C*-algebra (I,II)}.
\newblock K{\o}benhavns Universitet, Matematisk Institut, preprint no.5a/5b,
  1993.

\bibitem{Brown:85}
L.~G. Brown.
\newblock {Close hereditary C*-subalgebras and the structure of
  quasi-multipliers}.
\newblock MSRI preprint no. 11211-85, Purdue University, West Lafayette, 1985.

\bibitem{Brown:88}
L.~G. Brown.
\newblock {Semicontinuity and multipliers of C*-algebras}.
\newblock {\em Canad. J. Math.}, {\bf 40}, (1988).
\newblock 865-988.

\bibitem{Brown/Green/Rieffel:77}
L.~G. Brown, P.~Green, and M.~A. Rieffel.
\newblock {Stable isomorphism and strong Morita equivalence of C*-algebras}.
\newblock {\em Pacific J. Math.}, {\bf 71}, (1977).
\newblock 349-363.

\bibitem{Brown/Mingo/Shen:94}
L.~G. Brown, J.~A. Mingo, and N.~T. Shen.
\newblock {Quasi-multipliers and embeddings of Hilbert C*-bimodules}.
\newblock {\em Canad. J. Math.}, {\bf 46}, (1994).
\newblock 1150-1174.

\bibitem{Bui:92}
Huu~Hung Bui.
\newblock {\em {Morita equivalence of crossed products}}.
\newblock PhD thesis, University of New South Wales, Australia, 1992.

\bibitem{Bui:94/4}
Huu~Hung Bui.
\newblock {Induced representations twisted by cocycles}.
\newblock {\em Bull. Austral. Math. Soc.}, {\bf 50}, (1994).
\newblock 399-404.

\bibitem{Bui:94}
Huu~Hung Bui.
\newblock {Morita equivalence of twisted crossed products by coactions}.
\newblock {\em J. Functional Anal.}, {\bf 123}, (1994).
\newblock 59-98.

\bibitem{Bui:95/2}
Huu~Hung Bui.
\newblock {Full coactions on Hilbert C*-modules}.
\newblock {\em J. Austral. Math. Soc.}, {\bf 59}, (1995).
\newblock 409-420.

\bibitem{Bui:95/1}
Huu~Hung Bui.
\newblock {Morita equivalence of Twisted Crossed Products}.
\newblock {\em Proc. Amer. Math. Soc.}, {\bf 123}, (1995).
\newblock 2771-2776.

\bibitem{Bui:95/3}
Huu~Hung Bui.
\newblock {The crossed products of Hilbert C*-modules}.
\newblock report 95-185, Macquarie Univ., School of Mathematics, Physics,
  Computing and Electronics, New South Wales, Australia, 1995.

\bibitem{Bui:96/1}
Huu~Hung Bui.
\newblock {Full coactions on Hilbert C*-modules}.
\newblock preprint, Macquarie Univ., School of Mathematics, Physics, Computing
  and Electronics, New South Wales, Australia, 1996.

\bibitem{Bultheel:82}
A.~Bultheel.
\newblock {Inequalities in Hilbert modules of matrix-valued functions}.
\newblock {\em Proc. Amer. Math. Soc.}, {\bf 85}, (1982).
\newblock 369-372.

\bibitem{Bunke:92}
U.~Bunke.
\newblock {A K-theoretic relative index theorem}.
\newblock preprint MPI 92-55.
\newblock Max-Planck-Institut f\"ur Mathematik, Bonn, F.R.G., 1992.

\bibitem{Bunke:95}
U.~Bunke.
\newblock {A K-theoretic relative index theorem and Callias-type Dirac
  operators}.
\newblock {\em Math. Ann.}, {\bf 303}, (1995).
\newblock 241-279.

\bibitem{Bures:71}
D.~Bures.
\newblock {Abelian subalgebras of von Neumann algebras}.
\newblock {\em Memoirs Amer. Math. Soc.}, {\bf 110}, (1971).

\bibitem{Carey/Phillips:91}
A.~L. Carey and J.~Phillips.
\newblock {Algebras almost commuting with Clifford algebras in a ${\rm
  II}_\infty$ factor}.
\newblock {\em K-theory}, {\bf 5}, (1991).
\newblock 445-478.

\bibitem{Combes:84}
F.~Combes.
\newblock {Crossed products and Morita equivalence}.
\newblock {\em Proc. London Math. Soc.}, {\bf 49}, (1984).
\newblock 289-306.

\bibitem{Combes/Zettl:83}
F.~Combes and H.~H. Zettl.
\newblock {Order structures, traces and weights on Morita equivalent
  C*-algebras}.
\newblock {\em Math. Ann.}, {\bf 265}, (1983).
\newblock 67-81.

\bibitem{Connes:80}
A.~Connes.
\newblock {C*-alg\'ebres et g\'eom\'etrie diff\'erentielle}.
\newblock {\em C. R. Acad. Sci. Paris}, {\bf 290}, (1980).
\newblock 599-604.

\bibitem{Connes:80/2}
A.~Connes.
\newblock {Correspondances}.
\newblock Manuscriptus, 1980.

\bibitem{Connes:81}
A.~Connes.
\newblock {An analogue of the Thom isomorphism for crossed products of a
  C*-algebra by an action of $R$}.
\newblock {\em Adv. in Math.}, {\bf 39}, (1981).
\newblock 31-55.

\bibitem{Connes:82}
A.~Connes.
\newblock {A survey of foliations and operator algebras}.
\newblock {\em Proc. Symp. Pure Math. Amer. Math. Soc.}, {\bf 38}, (1982).
\newblock 521-628.

\bibitem{Connes:86}
A.~Connes.
\newblock {Non-commutative differential geometry}.
\newblock {\em Publ. I.H.E.S.}, {\bf 62}, (1985).
\newblock 41-144.

\bibitem{Connes:90/2}
A.~Connes.
\newblock {\em {G\'eom\'etrie non commutative}}.
\newblock Inter{E}ditions, Paris, 1990.

\bibitem{Connes:90}
A.~Connes.
\newblock {Introduction \`a la g\'eom\'etrie non-commutative}.
\newblock {\em Proc. Symp. Pure Math.}, {\bf 50}, (1990).
\newblock 91-118.

\bibitem{Connes:94}
A.~Connes.
\newblock {\em {Noncommutative geometry}}.
\newblock Academic Press, 1994.

\bibitem{Connes/Rieffel:87}
A.~Connes and M.~A. Rieffel.
\newblock {Yang--Mills for non-commutative two-tori}.
\newblock {\em Contemp. Math.}, {\bf 62}, (1987).
\newblock 237-266.

\bibitem{Connes/Skandalis:84}
A.~Connes and G.~Skandalis.
\newblock {The longitudinal index theorem for foliations}.
\newblock {\em Publ. RIMS, Kyoto Univ.}, {\bf 20}, (1984).
\newblock 1139-1183.

\bibitem{Cuntz:86}
J.~Cuntz.
\newblock {K-theory and C*-algebras}.
\newblock {\em Lect. Notes Math.}, {\bf 1046}, (1986).
\newblock 55-79.

\bibitem{Cuntz:93}
J.~Cuntz.
\newblock {A survey on some aspects of non-commutative geometry}.
\newblock {\em Jahresbericht der DMV}, {\bf 95}, (1993).
\newblock 60-84.

\bibitem{Cuntz/Higson:87}
J.~Cuntz and N.~Higson.
\newblock {Kuiper's theorem for Hilbert modules}.
\newblock {\em Contemp. Math.}, {\bf 62}, (1987).
\newblock 429-435.

\bibitem{Cuntz/Skandalis:86}
J.~Cuntz and G.~Skandalis.
\newblock {Mapping cones and exact sequences in KK-theory}.
\newblock {\em J. Operator Theory}, {\bf 15}, (1986).
\newblock 163-180.

\bibitem{Curto/Muhly/Williams:84}
R.~E. Curto, P.~S. Muhly, and D.~P. Williams.
\newblock {Crossed products of strongly Morita equivalent C*-algebras}.
\newblock {\em Proc. Amer. Math. Soc.}, {\bf 90}, (1984).
\newblock 528-530.

\bibitem{Dadarlat:88}
M.~D{\v{a}}d{\v{a}}rlat.
\newblock {On homomorphisms of matrix algebras of continuous functions}.
\newblock {\em Pacific J. Math.}, {\bf 132}, (1988).
\newblock 227-231.

\bibitem{Damir:94}
B.~Damir.
\newblock {Hilbert C*-modules over compact adjointable operators}.
\newblock In D.~Butkovi\'c et~al., editor, {\em {Functional analysis IV, Proc.
  of the postgraduate school and conference, Inter-University Center Dubrovnik,
  Croatia, Nov. 10-17, 1993}}, volume~{\bf 43} of {\em Var. Publ. Ser.} Aarhus
  University, Mat. Institut, Aarhus, Denmark, (1994).
\newblock 7-9.

\bibitem{DeBrabander:84}
M.~{De Brabander}.
\newblock {The classification of rational rotation algebras}.
\newblock {\em Arch. Math.}, {\bf 43}, (1984).
\newblock 79-83.

\bibitem{Delanghe:76}
R.~Delanghe.
\newblock {On Hilbert modules with reproducing kernel}.
\newblock {\em Lecture Notes Math.}, {\bf 561}, (1976).
\newblock 158-170.

\bibitem{Denizeau/Havet:93/2}
Y.~Denizeau and J.-F. Havet.
\newblock {Correspondances d'indice fini II: Indice d'une correspondance}.
\newblock preprint 93-04, Universit\'ed'Orl\'eans, France, 1993.

\bibitem{Denizeau/Havet:94}
Y.~Denizeau and J.-F. Havet.
\newblock {Correspondances d'indice fini I: Indice d'un vecteur}.
\newblock {\em J. Operator Theory}, {\bf 32}, (1994).
\newblock 111-156.

\bibitem{Dixmier/Douady:63}
J.~Dixmier and A.~Douady.
\newblock {Champs continus d'espaces hilbertiens et de C*-alg\`ebres}.
\newblock {\em Bull. Math. Soc. France}, {\bf 91}, (1963).
\newblock 227-283.

\bibitem{Douglas/Paulsen:89}
R.~G. Douglas and V.~I. Paulsen.
\newblock {\em {Hilbert modules over function algebras}}.
\newblock Pitman Res. Notes Math. v.217. Longman, New York, 1989.

\bibitem{Dupre:72}
M.~J. Dupr{\'e}.
\newblock {\em {The classification of Hilbert bundles}}.
\newblock PhD thesis, Univ. of Pennsylvania, U.S.A., 1972.

\bibitem{Dupre:74}
M.~J. Dupr{\'e}.
\newblock {Classifying Hilbert bundles, I}.
\newblock {\em J. Functional Analysis}, {\bf 15}, (1974).
\newblock 244-278.

\bibitem{Dupre:76}
M.~J. Dupr{\'e}.
\newblock {Classifying Hilbert bundles, II}.
\newblock {\em J. Functional Analysis}, {\bf 22}, (1976).
\newblock 295-322.

\bibitem{Dupre/Fillmore:81}
M.~J. Dupr{\'e} and P.~A. Fillmore.
\newblock {Triviality theorems for Hilbert modules}.
\newblock In {\em Topics in modern operator theory, 5th International
  conference on operator theory, Timi\c{s}oara and Herculane(Romania), June
  2-12, 1980}, Basel--Boston--Stuttgart, 1981. Birkh\"auser Verlag.
\newblock pp.71-79.

\bibitem{Dupre/Gillette:83}
M.~J. Dupr\'e and R.~M. Gillette.
\newblock {\em {Banach bundles, Banach modules and automorphisms of
  C*-algebras}}.
\newblock Research Notes in Mathematics v. 92, Adv. Publ. Program. Pitman,
  Boston - London - Melbourne, 1983.

\bibitem{Echterhoff:90}
S.~Echterhoff.
\newblock {On induced covariant systems}.
\newblock {\em Proc. Amer. Math. Soc.}, {\bf 108}, (1990).
\newblock 703-706.

\bibitem{Echterhoff:93}
S.~Echterhoff.
\newblock {Regularization of twisted covariant systems and crossed products
  with continuous trace}.
\newblock {\em J. Funct. Anal.}, {\bf 116}, (1993).
\newblock 277-313.

\bibitem{Echterhoff:94/3}
S.~Echterhoff.
\newblock {Duality of induction and restriction for abelian twisted covariant
  systems}.
\newblock {\em Math. Proc. Cambridge Philos. Soc.}, {\bf 119}, (1994).
\newblock 301-315.

\bibitem{Echterhoff:94}
S.~Echterhoff.
\newblock {Morita equivalent twisted actions and a new version of the
  Packer-Raeburn stabilization trick}.
\newblock {\em J. London Math. Soc. (2)}, {\bf 50}, (1994).
\newblock 170-186.

\bibitem{Echterhoff:94/2}
S.~Echterhoff.
\newblock {On transformation group C*-algebras with continuous trace}.
\newblock {\em Trans. Amer. Math. Soc.}, {\bf 343}, (1994).
\newblock 117-133.

\bibitem{Echterhoff:95}
S.~Echterhoff.
\newblock {Crossed products with continuous trace}.
\newblock {\em Memoirs Amer. Math. Soc.}, ???, 1995.
\newblock 136 pp.

\bibitem{Echterhoff/Raeburn:96/2}
S.~Echterhoff and I.~Raeburn.
\newblock {The stabilization trick for coactions}.
\newblock {\em J. Reine. Angew. Math.}, {\bf 470}, (1996).
\newblock 181-215.

\bibitem{Elliott:84}
G.~A. Elliott.
\newblock {On the K-theory of the C*-algebra generated by a projective
  representation of a torsion free discrete abelian group}.
\newblock In {\em {Operator Algebras and Group Representations, v. 1}}. Pitman,
  London, 1984.
\newblock 159-164.

\bibitem{Elliott/Qing:95}
G.~A. Elliott and Lin Qing.
\newblock {Cut-down method in the inductive limit decomposition of
  non-commutative tori}.
\newblock preprint, 1995.

\bibitem{Elliott/Saito/Wright:83}
G.~A. Elliott, K.~Sait{\^o}, and J.~D.~M. Wright.
\newblock {Embedding AW*-algebras as double commutants in type I algebras}.
\newblock {\em J. London Math. Soc.}, {\bf 28}, (1983).
\newblock 376-384.

\bibitem{Ellis/Gohberg/Lay:95}
R.~L. Ellis, I~Gohberg, and D.~C. Lay.
\newblock {Infinite analogues of block Toeplitz matrices and related orthogonal
  functions}.
\newblock {\em Integral Equations Oper. Theory}, {\bf 22}, (1995).
\newblock 375-419.

\bibitem{Exel:93}
R.~Exel.
\newblock {A Fredholm operator approach to Morita equivalence}.
\newblock {\em K-Theory}, {\bf 7}, (1992).
\newblock 285-308.

\bibitem{Exel:94}
R.~Exel.
\newblock {Twisted partial actions -- a classification of stable C*-algebraic
  bundles}.
\newblock preprint, Univ. de S\~ao Paulo, Brazil, 1994.

\bibitem{Fell:77}
J.~Fell.
\newblock {\em {Induced representations and Banach $\ast$-algebraic bundles}}.
\newblock Lecture Notes Math. v. {\bf 582}. Springer--Verlag,
  Berlin--Heidelberg--New York, (1977).

\bibitem{Fell/Doran:88}
J.~M.~G. Fell and R.~S. Doran.
\newblock {\em {Representation of $\ast$-algebras, locally compact groups, and
  Banach $\ast$-algebraic bundles, I+II}}.
\newblock Pure and Applied Mathematics, v. {\bf 125} and {\bf 126}. Academic
  Press, Inc., Boston, Mass., 1988.

\bibitem{Filippov:87}
O.~G. Filippov.
\newblock {On C*-algebras A over which the Hilbert module $l_2(A)$ is self-dual
  (russ./engl.)}.
\newblock {\em Vestn. Mosk. Univ., Ser. I: Mat.-Mekh.}, no. 4, (1987).
\newblock 74-76 / {\it Moscow Univ. Math. Bull.} {\bf 42}(1987), no. 4, 87-90.

\bibitem{Filippov:90/1}
O.~G. Filippov.
\newblock {On the reflexivity of objects of certain concrete categories
  (russ./engl.)}.
\newblock {\em Vestn. Mosk. Univ., Ser. I: Mat.-Mekh.}, no. 1, (1990).
\newblock 93-95 / {\it Moscow Univ. Math. Bull.} {\bf 45}(1990), no. 1, 53-54.

\bibitem{Filippov:90/2}
O.~G. Filippov.
\newblock {The reduction of operators with almost periodic symbols to operators
  over C*-algebras on sections of associated bundles over a torus}.
\newblock {\em Ann. Global Anal. Geom.}, {\bf 8}, (1990).
\newblock 113-126.

\bibitem{Fox/Haskell/Raeburn:89}
J.~Fox, P.~Haskell, and I.~Raeburn.
\newblock {Kasparov products, KK-equivalence and proper actions of connected
  Lie groups}.
\newblock {\em J. Oper. Theory}, {\bf 22}, (1989).
\newblock 3-29.

\bibitem{Frank:85}
M.~Frank.
\newblock {A set of maps from $K$ to $End_A(l_2(A))$ isomorphic to
  $End_{A(K)}(l_2(A(K)))$. Applications.}
\newblock {\em Annals Global Anal. Geom.}, {\bf 3}, (1985).
\newblock 155-171.

\bibitem{Frank:88}
M.~Frank.
\newblock {\em {Beitr\"age zur Entwicklung und systematischen Darstellung der
  Theorie der Hilbert-C*-Moduln}}.
\newblock PhD thesis, Karl--Marx--Universit\"at Leipzig, G.D.R., 1988.

\bibitem{Frank:89/2}
M.~Frank.
\newblock {Central direct integral decomposition of von Neumann algebras and
  some operator algebras on self-dual Hilbert W*-modules over commutative
  W*-algebras}.
\newblock preprint no. 13, KMU--CLG, Leipzig, G.D.R., 1989.

\bibitem{Frank:89/1}
M.~Frank.
\newblock {Elements of Tomita-Takesaki theory for embedable AW*-algebras}.
\newblock {\em Annals Global Anal. Geom.}, {\bf 7}, (1989).
\newblock 115-131.

\bibitem{Frank:90/2}
M.~Frank.
\newblock {One-parameter groups arising from some real subspaces of self-dual
  Hilbert W*-modules}.
\newblock {\em Math. Nachr.}, {\bf 145}, (1990).
\newblock 169-185.

\bibitem{Frank:90/1}
M.~Frank.
\newblock {Self-duality and C*-reflexivity of Hilbert C*-modules}.
\newblock {\em Zeitschr. Anal. Anw.}, {\bf 9}, (1990).
\newblock 165-176.

\bibitem{Frank:90/3}
M.~Frank.
\newblock {Von Neumann representations on self-dual Hilbert W*-modules}.
\newblock {\em Math. Nachr.}, {\bf 145}, (1990).
\newblock 187-199.

\bibitem{Frank:92/6}
M.~Frank.
\newblock {Direct integrals and Hilbert W*-modules (russ.)}.
\newblock In {\em Problems in algebra, geometry and discrete mathematics, eds.:
  O. B. Lupanov, A. I. Kostrikin, Moscow State University, Dept. Mech. Math.,
  Moscow, Russia}, (1992).
\newblock 162-177.

\bibitem{Frank:93/5}
M.~Frank.
\newblock {Geometrical aspects of Hilbert C*-modules}.
\newblock K{\o}benhavns Universitet, Matematisk Institut, preprint 22/1993,
  Copenhagen, Denmark, 1993.

\bibitem{Frank:93/1}
M.~Frank.
\newblock {Normal operator-valued weights of finite index}.
\newblock preprint no. 8/93, NTZ, Univ. Leipzig, F.R.G., 1993.

\bibitem{Frank:95/4}
M.~Frank.
\newblock {Hilbert C*-modules over monotone complete C*-algebras}.
\newblock {\em Math. Nachr.}, {\bf 175}, (1995).
\newblock 61-83.

\bibitem{Frank:95/3}
M.~Frank.
\newblock {Isomorphisms of Hilbert C*-modules and $*$-isomorphisms of related
  operator C*-algebras}.
\newblock preprint no. 14/95, ZHS-NTZ, Univ. Leipzig, F.R.G. / accepted by
  Math. Scand., (to appear in fall 1997), 1995.

\bibitem{Frank:96/3}
M.~Frank.
\newblock {On the Hahn-Banach theorem for Hilbert C*-modules}.
\newblock preprint 11/1996, ZHS-NTZ, Univ. Leipzig, F.R.G., 1996.

\bibitem{Frank/Manuilov:95}
M.~Frank and V.~M. Manuilov.
\newblock {Diagonalizing ''compact'' operators on Hilbert W*-modules}.
\newblock {\em Zeitschr. Anal. Anwendungen}, {\bf 14}, (1995).
\newblock 33-41.

\bibitem{Frank/Troitsky:95}
M.~Frank and E.~V. Troitsky.
\newblock {Lefschetz numbers and geometry of operators in W*-modules}.
\newblock preprint 4/95, ZHS-NTZ, Univ. Leipzig, F.R.G. / submitted to {\it
  Funkts.~Anal.~i Prilozh.}, Moscow, 1995.

\bibitem{Frydryszak/Jakobczyk:88}
A.~Frydryszak and L.~Jakobczyk.
\newblock {Generalized Gelfand-Naimark-Segal construction for supersymmetric
  quantum mechanics}.
\newblock {\em Lett. Math. Phys.}, {\bf 16}, (1988).
\newblock no. 2, 101-107.

\bibitem{Fuge:95}
H.~Fuge.
\newblock {Einige Aussagen \"uber beschr\"ankte Operatoren in Hilbertmoduln}.
\newblock Master's thesis, Univ. Leipzig, 1995.

\bibitem{Gekhtman:84}
M.~I. Gekhtman.
\newblock {Hilbert modules and pseudo-Hilbert spaces, (russian)}.
\newblock In {\em {Spectral theory of operators and infinite analysis, Collect.
  sci. works}}. Kiew, USSR, 1984.
\newblock 57-65.

\bibitem{Gerecke/Lorenz:95}
U.~Gerecke and J.~Lorenz.
\newblock {Grundlegende Aussagen \"uber nichtnegativ hermitesche Ma{\ss}e,
  Ma{\ss}e mit orthogonalen Werten sowie projektorwertige Ma{\ss}e}.
\newblock Master's thesis, Universit\"at Leipzig, Leipzig, F.R.G., 1995.
\newblock 367 pp.

\bibitem{Ghez/Lima/Roberts:85}
P.~Ghez, R.~Lima, and J.~E. Roberts.
\newblock {W*-categories}.
\newblock {\em Pacific J. Math.}, {\bf 120}, (1985).
\newblock 79-109.

\bibitem{Giordano:88}
T.~Giordano.
\newblock {A classification of approximately finite real C*-algebras}.
\newblock {\em J. Reine Angew. Math.}, {\bf 385}, (1988).
\newblock 161-194.

\bibitem{Giordano/Handelman:89}
T.~Giordano and D.~E. Handelman.
\newblock {Real AF C*-algebras with $K_o$ of small rank}.
\newblock {\em Can. J. Math.}, {\bf 41}, (1989).
\newblock 786-807.

\bibitem{Goldstine/Horwitz:66}
H.~H. Goldstine and L.~P. Horwitz.
\newblock {Hilbert space with non-associative scalars}.
\newblock {\em Math. Ann.}, {\bf 164}, (1966).
\newblock 291-316.

\bibitem{Green:78}
P.~Green.
\newblock {The local structure of twisted covariance algebras}.
\newblock {\em Acta Math.}, {\bf 140}, (1978).
\newblock 191-250.

\bibitem{Greene:74}
W.~A. Greene.
\newblock {Ambrose modules}.
\newblock {\em Mem. Amer. Math. Soc.}, {\bf 148}, (1974).
\newblock 109-134.

\bibitem{Haagerup:89}
U.~Haagerup.
\newblock {The injective factors of type $III_\lambda$, $ 0< \lambda <1$}.
\newblock {\em Pacific J. Math.}, {\bf 137}, (1989).
\newblock 265-310.

\bibitem{Halpern:85}
H.~Halpern.
\newblock {One parameter automorphism groups of generalized KH-algebras}.
\newblock preprint, Univ. of Cincinnati, Cincinnati, USA, 1985.

\bibitem{Hamana:79/2}
M.~Hamana.
\newblock {Injective envelopes of C*-algebras}.
\newblock {\em J. Math. Soc. Japan}, {\bf 31}, (1979).
\newblock 181-197.

\bibitem{Hamana:79/1}
M.~Hamana.
\newblock {Injective envelopes of operator systems}.
\newblock {\em Publ. Res. Inst. Math. Sci. Kyoto}, {\bf 15}, (1979).
\newblock 773-785.

\bibitem{Hamana:81}
M.~Hamana.
\newblock {Regular embeddings of C*-algebras in monotone complete C*-algebras}.
\newblock {\em J. Math. Soc. Japan}, {\bf 33}, (1981).
\newblock 159-183.

\bibitem{Hamana:82/2}
M.~Hamana.
\newblock {Tensor products for monotone complete C*-algebras}.
\newblock {\em Japan. J. Math.}, {\bf 8}, (1982).
\newblock 259-283.

\bibitem{Hamana:82/3}
M.~Hamana.
\newblock {Tensor products for monotone complete C*-algebras, II}.
\newblock {\em Japan. J. Math.}, {\bf 8}, (1982).
\newblock 285-295.

\bibitem{Hamana:92}
M.~Hamana.
\newblock {Modules over monotone complete C*-algebras}.
\newblock {\em Internat. J. Math.}, {\bf 3}, (1992).
\newblock 185-204.

\bibitem{Havet:88}
J.-F. Havet.
\newblock {Calcul fonctonnel continu dans les modules hilbertiens autoduaux}.
\newblock preprint 1988, Or\'eans, France.

\bibitem{Havet:90}
J.-F. Havet.
\newblock {Esp\'erance conditionelle minimale}.
\newblock {\em J. Oper. Theory}, {\bf 24}, (1990).
\newblock 33-55.

\bibitem{Hegerfeldt:85}
Gerhard~C. Hegerfeldt.
\newblock {Inequalities of Schwarz and H\"older type for random operators}.
\newblock {\em J. Math. Phys.}, {\bf 26}, (1985).
\newblock 1576-1577.

\bibitem{Hennings:89}
M.~A. Hennings.
\newblock {Kasparov's technical lemma for b*-algebras}.
\newblock {\em Math. Proc. Cambridge Philos. Soc.}, {\bf 105}, (1989).
\newblock 537-545.

\bibitem{Heo:96}
Jaeseong Heo.
\newblock {Completely multi-positive linear maps and representations on Hilbert
  C*-modules}.
\newblock preprint, Seoul National University, Korea, 1996.

\bibitem{Higson:88}
N.~Higson.
\newblock {Algebraic K-theory of stable C*-algebras}.
\newblock {\em Adv. Math.}, {\bf 67}(no. 1), (1988).
\newblock 140pp.

\bibitem{Higson:90}
N.~Higson.
\newblock {A primer on KK--theory}.
\newblock {\em Proc. Symp. Pure Math.}, {\bf 51-1}, (1990).
\newblock 239-285.

\bibitem{Hilsum:85}
M.~Hilsum.
\newblock {Signature operator on Lipschitz manifolds and unbounded Kasparov
  bimodules}.
\newblock In H.~Araki, C.~C. Moore, \c{S}. Stratila, and D.~Voiculescu,
  editors, {\em Lecture Notes Math. v. {\bf 1132}, Operator algebras and their
  connections with topology and ergodic theory}. Springer--Verlag, Berlin,
  1985.
\newblock pp. 254-288.

\bibitem{Hilsum:89}
M.~Hilsum.
\newblock {Fonctorialit\'e en K-th\'eorie bivariante pour les vari\'et\'es
  lipschitziennes}.
\newblock {\em K-theory}, {\bf 3}, (1989).
\newblock 401-440.

\bibitem{Hilsum/Skandalis:83}
M.~Hilsum and G.~Skandalis.
\newblock {Stabilit\'e des C*-alg\'ebres de feuilletages}.
\newblock {\em Ann. Inst. Fourier, Grenoble}, {\bf 33}, (1983).
\newblock 201-208.

\bibitem{Hilsum/Skandalis:92}
M.~Hilsum and G.~Skandalis.
\newblock {Invariance par homotopie de la signature a coefficients dans un
  fibre presque plat}.
\newblock {\em J. Reine Angew. Math.}, {\bf 423}, (1992).
\newblock 73-99.

\bibitem{Hofmann:72}
K.~H. Hofmann.
\newblock {Representations of algebras by continuous sections}.
\newblock {\em Bull. Amer. Math. Soc.}, {\bf 78}, (1972).
\newblock 291-373.

\bibitem{Hofmann:74}
K.~H. Hofmann.
\newblock {Erratum : Representations of algebras by continuous sections}.
\newblock {\em Mem. Amer. Math. Soc.}, {\bf 148}, (1974).
\newblock 177-182.

\bibitem{Irmatov:88}
A.~A. Irmatov.
\newblock {On a topology in the space of Fredholm operators (russ.)}.
\newblock In {\em {Selected questions of algebra, geometry and discrete
  mathematics, Moscow, Russia, 40-43}}, 1988.

\bibitem{Irmatov:89}
A.~A. Irmatov.
\newblock {On a new topology in the space of Fredholm operators}.
\newblock {\em Ann. Global Anal. Geom.}, {\bf 7}, (1989).
\newblock 93-106.

\bibitem{Irmatov:90}
A.~A. Irmatov.
\newblock {The topology of the space of Fredholm operators and invariants of
  non-linear Fredholm maps (russ.)}.
\newblock {\em Uspekhi Mat. Nauk}, {\bf 45}(1), (1990).
\newblock 173-174.

\bibitem{Irmatov/Mishchenko:90}
A.~A. Irmatov and A.~S. Mi\v{s}\v{c}enko.
\newblock {Infinitesimal Fredholm structures on infinite-dimensional
  manifolds}.
\newblock In W.~B. Arveson, A.~S. Mi\v{s}\v{c}enko, M.~Putinar, M.~A. Rieffel,
  and S.~\v{S}tratila, editors, {\em Operator algebras and Topology, Proc. of
  the OATE 2 Conf., Romania 1989}, Pitman Research Notes in Mathematics Series
  v. 270. Longman Scientific \& Technical, New York, 1990.
\newblock 45-81.

\bibitem{Istratescu:87}
V.~I. Istr\v{a}\c{t}escu.
\newblock {\em {Inner product structures, Theory and Applications}}.
\newblock Mathematics and Its Applications v.25. D. Reidel Publishing Company,
  Dordrecht - Boston - Lancaster - Tokyo, 1987.

\bibitem{Itoh:80}
S.~Itoh.
\newblock {A note on dilations in modules over C*-algebras}.
\newblock {\em J. London Math. Soc.}, {\bf 22}, (1980).
\newblock 117-126.

\bibitem{Itoh:90}
S.~Itoh.
\newblock {Reproducing kernels in modules over C*-algebras and their
  applications}.
\newblock {\em Bull. Kyushu Inst. Technol. Math. Natur. Sci.}, {\bf 37},
  (1990).
\newblock 1-20.

\bibitem{Je/Yang:86}
Hai-Gon Je and Young-Oh Yang.
\newblock {On the spatial numerical ranges and Hermitian operators}.
\newblock {\em Univ. Ulsan Rep. Natur. Sci. Eng.}, {\bf 17}, (1986).
\newblock no. 2, 209-214.

\bibitem{Jensen/Thomsen:91}
K.~K. Jensen and K.~Thomsen.
\newblock {\em {Elements of KK-theory}}.
\newblock Mathematics: Theory and Applications. Birkh\"auser, Boston, Mass.,
  1991.

\bibitem{Jeong:93}
Ja~A Jeong.
\newblock {Full hereditary C*-subalgebras of crossed products}.
\newblock {\em Bull. Korean Math. Soc.}, {\bf 30}, (1993).
\newblock 193-199.

\bibitem{Ji:94}
G.~Ji.
\newblock {Generalized Cowen-Douglas operators over Hilbert C*-modules}.
\newblock {\em Integr. Equat. Oper. Th.}, {\bf 20}, (1994).
\newblock 395-409.

\bibitem{Jolissaint:91}
P.~Jolissaint.
\newblock {Indice d'esp\'erances conditionelles et alg\`ebres de von Neumann
  finies}.
\newblock {\em Math. Scand.}, {\bf 68}, (1991).
\newblock 221-246.

\bibitem{Kaftal:91}
V.~Kaftal.
\newblock {Type decomposition for von Neumann algebra embeddings}.
\newblock {\em J. Funct. Anal.}, {\bf 98}, (1991).
\newblock 169-193.

\bibitem{Kajiwara:87}
T.~Kajiwara.
\newblock {Remarks on strongly Morita equivalent C*-crossed products}.
\newblock {\em Math. Jap.}, {\bf 32}, (1987).
\newblock 257-260.

\bibitem{Kajiwara/Watatani:95}
T.~Kajiwara and Y.~Watatani.
\newblock {Jones index theory by Hilbert C*-bimodules and K-theory}.
\newblock preprint, Dept.~Environm.~Math.~Sci., Okayama Univ., Japan, 1995.

\bibitem{Kakihara:80}
Y.~Kakihara.
\newblock {Hilbert B(H)-modules with applications, II}.
\newblock {\em Res. Rep. Inst. Inform. Sci. Tech.}, {\bf 6}, (1980).
\newblock 33-45.

\bibitem{Kakihara:82}
Y.~Kakihara.
\newblock {Hilbert B(H)-modules with applications, III}.
\newblock {\em Res. Activ. Fac. Sci. and Eng., Tokyo Denki Univ.}, {\bf 4},
  (1982).
\newblock 11-32.

\bibitem{Kakihara:83/1}
Y.~Kakihara.
\newblock {Hilbert B(H)-modules with applications, IV}.
\newblock {\em Res. Activ. Fac. Sci. and Eng., Tokyo Denki Univ.}, {\bf 5},
  (1983).
\newblock 27-30.

\bibitem{Kakihara:83/2}
Y.~Kakihara.
\newblock {On a Hilbert module over an operator algebra and its application to
  harmonic analysis}.
\newblock {\em Kodai Math. J.}, {\bf 6}, (1983).
\newblock 289-300.

\bibitem{Kakihara:84}
Y.~Kakihara.
\newblock {Hilbert B(H)-modules with applications, V}.
\newblock {\em Res. Activ. Fac. Sci. and Eng., Tokyo Denki Univ.}, {\bf 6},
  (1984).
\newblock 131-132.

\bibitem{Kakihara:85}
Y.~Kakihara.
\newblock {A note on harmonizable and V-bounded processes}.
\newblock {\em J. Multivar. Anal.}, {\bf 16}, (1985).
\newblock 140-156.

\bibitem{Kakihara/Teresaki:79}
Y.~Kakihara and T.~Teresaki.
\newblock {Hilbert B(H)-modules with applications, I}.
\newblock {\em Res. Rep. Inst. Inf. Sci. Tech.}, {\bf 5}, (1979).
\newblock 23-32.

\bibitem{Kaliszewski:94}
S.~P. Kaliszewski.
\newblock {\em {Morita equivalence methods for twisted C*-dynamical systems}}.
\newblock PhD thesis, Dartmouth College, New Hampshire, USA, 1994.

\bibitem{Kaliszewski:95}
S.~P. Kaliszewski.
\newblock {A note on Morita equivalence of twisted C*-dynamical systems}.
\newblock {\em Proc. Amer. Math. Soc.}, 123, (1995).
\newblock 1737-1740.

\bibitem{Kaliszewski/Quigg:95}
S.~P. Kaliszewski and J.~Quigg.
\newblock {Imprimitivity for C*-coactions of non-amenable groups}.
\newblock preprint, funct-an@babbage.sissa.it, 9602003, 1995.

\bibitem{Kaliszewski/Quigg/Raeburn:95}
S.~P. Kaliszewski, J.~Quigg, and I.~Raeburn.
\newblock {Duality of restriction and induction for C*-coactions}.
\newblock preprint, funct-an@babbage.sissa.it, 9602004, 1995.

\bibitem{Kaminker/Miller:85}
J.~Kaminker and J.~G. Miller.
\newblock {Homotopy invariance of the analytic index of signature operators
  over C*-algebras}.
\newblock {\em J. Oper. Theory}, {\bf 14}, (1985).
\newblock 113-127.

\bibitem{Kandelaki:86}
T.~K. Kandelaki.
\newblock {Category of homomorphisms into a generalized Calkin algebra and
  projective modules over the commutant (russ.)}.
\newblock {\em Soobshch. Akad. Nauk Gruzin. SSR}, {\bf 122}, (1986).
\newblock 253-255.

\bibitem{Kaplansky:53}
I.~Kaplansky.
\newblock {Modules over operator algebras}.
\newblock {\em Amer. J. Math.}, {\bf 75}, (1953).
\newblock 839-858.

\bibitem{Kasimov:81}
V.~A. Kasimov.
\newblock {Hilbert structures on modules over C*-algebras (russ.)}.
\newblock {\em Akad. Nauk Azerba{\H{\i}}dshan. SSR Dokl.}, {\bf 37}, (1981).
\newblock 3-5.

\bibitem{Kasimov:82/2}
V.~A. Kasimov.
\newblock {A property of Hilbert modules and Fredholm operators over
  C*-algebras (russ./engl.)}.
\newblock {\em Akad. Nauk Azerba{\H{\i}}dshan. SSR Dokl.}, {\bf 38}, (1982).
\newblock 10-14 / {\it Amer. Math. Soc. Transl., Series 2} {\bf 136}(1987),
  143-147.

\bibitem{Kasimov:82/1}
V.~A. Kasimov.
\newblock {Homotopy properties of the general linear group of the Hilbert
  module $l_2(A)$ (russ./engl.)}.
\newblock {\em Mat. Sbornik}, {\bf 119}, (1982).
\newblock 376-386 / {\it Math. USSR - Sb.} {\bf 47}(1984), 365-376.

\bibitem{Kasimov:89/2}
V.~A. Kasimov.
\newblock {Homotopy properties of Hilbert modules (russ.)}.
\newblock {\em Studies in algebra and topology. Themat. Collect. Sci. Works,
  Baku}, 52-55, 1989.

\bibitem{Kasimov:89/1}
V.~A. Kasimov.
\newblock {Homotopy triviality of the group $GL^*(l_2(A))$ (russ.)}.
\newblock {\em Studies in algebra and topology. Themat. Collect. Sci. Works,
  Baku}, 46-51, 1989.

\bibitem{Kasparov:75}
G.~G. Kasparov.
\newblock {Topological invariants of elliptic operators. I: K-homology
  (russ./engl.)}.
\newblock {\em Izv. Akad. Nauk SSSR, Ser. Mat.}, {\bf 39}, (1975).
\newblock 796-838 / {\it Math. USSR - Izv.} {\bf 9}(1975), 751-792.

\bibitem{Kasparov:80}
G.~G. Kasparov.
\newblock {Hilbert C*-modules: theorems of Stinespring and Voiculescu}.
\newblock {\em J. Operator Theory}, {\bf 4}, (1980).
\newblock 133-150.

\bibitem{Kasparov:81}
G.~G. Kasparov.
\newblock {Operator K-theory and extensions of C*-algebras (russ./engl.)}.
\newblock {\em Izv. Akad. Nauk SSSR, Ser. Mat.}, {\bf 44}, (1980).
\newblock 571-630 / {\it Math. USSR - Izv.} {\bf 16}(1981), no. 3.

\bibitem{Kasparov:82}
G.~G. Kasparov.
\newblock {Conspectus}.
\newblock unpublished, 1982.

\bibitem{Kasparov:83}
G.~G. Kasparov.
\newblock {Operator K-theory and its applications: Elliptic Operators, Group
  Representations, Higher Signatures, C*-extensions}.
\newblock {\em Proc. Int. Congress of Math.}, 1983.
\newblock 987-1000.

\bibitem{Kasparov:88}
G.~G. Kasparov.
\newblock {Equivariant KK-theory and the Novikov conjecture}.
\newblock {\em Invent. Math.}, {\bf 91}, (1988).
\newblock 147-201.

\bibitem{Kasparov:93}
G.~G. Kasparov.
\newblock {Novikov's conjecture on higher signatures: the operator K-theory
  approach}.
\newblock {\em Contemp. Math.}, {\bf 145}, (1993).
\newblock 79-99.

\bibitem{Kasparov:85}
G.~G. Kasparov.
\newblock {Operator K-theory and its applications (russ./engl.)}.
\newblock {\em Itogi Nauki i Tekhn., Ser. Sovrem. Probl. Mat.}, {\bf 27},
  VINITI, Moscow, 1985.
\newblock 3-31 / {\it J. Soviet Math.} {\bf 37}(1987), 1373-1396.

\bibitem{Kasparov/Skandalis:90}
G.~G. Kasparov and G.~Skandalis.
\newblock {Groupes agissant sur des immeubles de Bruhat--Tits, K-th\'eorie
  operationelle et conjecture de Novikov}.
\newblock {\em C. R. Acad. Sci. Paris, Ser. I}, {\bf 310}, (1990).
\newblock 171-174.

\bibitem{Kasparov/Skandalis:91}
G.~G. Kasparov and G.~Skandalis.
\newblock {Groups acting on buildings, operator K-theory, and Novikov's
  conjecture}.
\newblock {\em K-theory}, {\bf 4}, (1991).
\newblock 303-338.

\bibitem{Khimshiashvili:92}
G.~N. Khimshiashvili.
\newblock {On homotopic structure of invertible singular operators}.
\newblock In G.~S. Chogoshvili, editor, {\em {Generalized Homologies and
  Homotopies. Collection of works in homology theory, 5. Tbilisi: Metsniereba
  (ISBN 5-520-00649-0)}}, volume~{\bf 97} of {\em {\it Tr. Tbilis. Mat. Inst.
  A. M. Razmadze}}, (1992).
\newblock 78-91.

\bibitem{Kim/Yang:84}
Kyong~Soo Kim and Youngoh Yang.
\newblock {On the numerical range for nonlinear operators}.
\newblock {\em Bull. Korean Math. Soc.}, {\bf 21}, (1984).
\newblock 119-126.

\bibitem{Kirchberg:91}
E.~Kirchberg.
\newblock {Commutants of unitaries in UHF-algebras and functorial properties of
  exactness}.
\newblock preprint, Heidelberg, F.~R.~Germany, 1991.

\bibitem{Kodaka:89/2}
K.~Kodaka.
\newblock {Automorphisms, diffeomorphisms and strong Morita equivalence of
  irrational rotation C*-algebras}.
\newblock {\em Tokyo J. Math.}, {\bf 12}, (1989).
\newblock 415-427.

\bibitem{Kodaka:89/1}
K.~Kodaka.
\newblock {Automorphisms of unital C*-algebras are strongly Morita equivalent
  to irrational rotation algebras}.
\newblock {\em Tokyo J. Math.}, {\bf 12}, (1989).
\newblock 175-179.

\bibitem{Kokschal:95}
A.~Kokschal.
\newblock {Einige Grundz\"uge der Theorie der Hilbertmoduln}.
\newblock Master's thesis, Univ. Leipzig, 1995.
\newblock 379 pp.

\bibitem{Kumijan:88}
A.~Kumijan.
\newblock {On equivariant sheaf cohomology and elementary C*-bundles}.
\newblock {\em J. Oper. Theory}, {\bf 20}, (1988).
\newblock 207-240.

\bibitem{Kurmakaev/Shkarin:90}
E.~Sh. Kurmakaeva and S.~A. Shkarin.
\newblock {On projectivity of some modules over polynormed algebras of
  continuous functions (russ.)}.
\newblock {\em Vestn. Mosk. Univ.}, no. 5, 1990.
\newblock 66-68.

\bibitem{Lance:94}
E.~C. Lance.
\newblock {Unitary operators on Hilbert C*-modules}.
\newblock {\em Bull. London Math. Soc.}, {\bf 26}, (1994).
\newblock 363-366.

\bibitem{Lance:95}
E.~C. Lance.
\newblock {\em {Hilbert C*-modules - a toolkit for operator algebraists}}.
\newblock London Mathematical Society Lecture Note Series {\bf 210}. Cambridge
  University Press, Cambridge, England, 1995.

\bibitem{Landsman:94}
N.~P. Landsman.
\newblock {Quantisierung der 'Moment map' durch Hilbert-C*-Moduln und
  Anwendungen in der algebraischen Quantenfeldtheorie}.
\newblock Verhandlungen der Deutschen Physikalischen Gesellschaft, Heft 4/1994,
  Bericht \"uber die 58. Physikerjahrestagung der Fachgremien, Kurzfassung des
  Vortrages MP 3A.12, p. 556, 1994.

\bibitem{Landsman:95}
N.~P. Landsman.
\newblock {Rieffel induction as generalized quantum Marsden-Weinstein
  reduction}.
\newblock {\em J. Geom. Phys.}, {\bf 15}, (1995).
\newblock 285-319.

\bibitem{Lee/Kim:82}
Sanghun Lee and Yongchan Kim.
\newblock {Double B-centralizers of pre--Hilbert B-modules}.
\newblock {\em Kyungpook Math. J.}, {\bf 22}, (1982).
\newblock 303-307.

\bibitem{Lesch:88}
M.~Lesch.
\newblock {\em {Die K-Theorie der C*-Algebra der Toeplitzoperatoren auf den
  Lie-Sph\"aren}}.
\newblock PhD thesis, Universit\"at Marburg/Lahn, FRG, 1988.

\bibitem{Lin:91/3}
H.~Lin.
\newblock {Bounded module maps and pure completely positive maps}.
\newblock {\em J. Operator Theory}, {\bf 26}, (1991).
\newblock 121-138.

\bibitem{Lin:91}
H.~Lin.
\newblock {Generalized Weyl - von Neumann theorems}.
\newblock {\em Int. J. Math.}, {\bf 2}, (1991).
\newblock 725-739.

\bibitem{Lin:91/2}
H.~Lin.
\newblock {Hilbert C*-modules and their bounded module maps}.
\newblock {\em Science in China (Series A)}, {\bf 34}(4), (1991).
\newblock 2-13.

\bibitem{Lin:92}
H.~Lin.
\newblock {Injective Hilbert C*-modules}.
\newblock {\em Pacific J. Math.}, {\bf 154}, (1992).
\newblock 131-164.

\bibitem{Lin:93}
H.~Lin.
\newblock {Extensions of C*-algebras with real rank zero}.
\newblock {\em Internat. J. Math.}, {\bf 4}, (1993).
\newblock 231-252.

\bibitem{Lin:93/2}
H.~Lin.
\newblock {Extensions of multipliers and injective Hilbert modules}.
\newblock {\em Chin. Ann. Math., Ser. B}, {\bf 14}, (1993).
\newblock 387-396.

\bibitem{Lin:94}
H.~Lin.
\newblock {The generalized Weyl - von Neumann theorem and C*-algebra
  extensions}.
\newblock In R.~Curto and P.~E.~T. J{\o}rgensen, editors, {\em Algebraic
  Methods in Operator Theory}. Birkh\"auser, Boston - Basel - Berlin, 1994.

\bibitem{Loginov/Shulman:93}
A.~I. Loginov and V.~S. Shulman.
\newblock {Vector-valued duality for modules over Banach algebras
  (russ./engl.)}.
\newblock {\em Izv. Akad. Nauk Rossii, Ser. Mat.}, {\bf 57}, (1993).
\newblock no. 4, 3-35 / {\it Russ. Acad. Sci. Izv., Math.} {\bf 43}(1994),
  no.~1, 1-29.

\bibitem{Lott/Connes:92}
J.~Lott and A.~Connes.
\newblock {The metric aspects of noncommutative geometry}.
\newblock Max-Planck-Institut f\"ur Mathematik, Bonn, FRG, preprint MPI 92-14,
  1992.

\bibitem{Loynes:65}
R.~M. Loynes.
\newblock {Linear operators in VH-spaces}.
\newblock {\em Trans. Amer. Math. Soc.}, {\bf 166}, (1965).
\newblock 167-180.

\bibitem{Loynes:65/3}
R.~M. Loynes.
\newblock {On a generalization of second-order stationary processes}.
\newblock {\em Proc. London Math. Soc.}, {\bf 15}, (1965).
\newblock 385-398.

\bibitem{Loynes:65/2}
R.~M. Loynes.
\newblock {On generalized positive definite functions}.
\newblock {\em Proc. London Math. Soc.}, {\bf 15}, (1965).
\newblock 373-384.

\bibitem{Lu:95}
Yun-Gang Lu.
\newblock {Passage from quantum systems with continuous spectrum to quantum
  Poisson processes on Hilbert modules}.
\newblock {\em J. Math. Phys.}, {\bf 36}, (1995).
\newblock 142-176.

\bibitem{Lutz:95}
F.~H. Lutz.
\newblock {Beispiele nichtkommutativer Geometrien}.
\newblock Diplomarbeit, Eberhard-Karls-Universit\"at T\"ubingen, F.R.G.,
  M\"arz, 1995.

\bibitem{Magajna:94}
B.~Magajna.
\newblock {Hilbert modules and completely bounded operators}.
\newblock preprint, University of Ljubljana, Ljubljana, Slovenia, 1994.

\bibitem{Magajna:95/2}
B.~Magajna.
\newblock {Hilbert modules and tensor products of operator spaces}.
\newblock preprint, University of Ljubljana, Ljubljana, Slovenia, 1995.

\bibitem{Magajna:95/1}
B.~Magajna.
\newblock {Tensor products over abelian W*-algebras}.
\newblock preprint, University of Ljubljana, Ljubljana, Slovenia, 1995.

\bibitem{Mallios:85}
A.~Mallios.
\newblock {Hermitian K-theory over topological $\ast$-algebras}.
\newblock {\em J. Math. Anal. Appl.}, {\bf 106}, (1985).
\newblock 454-539.

\bibitem{Mansfield:91}
K.~Mansfield.
\newblock {Induced representations of crossed products by coactions}.
\newblock {\em J. Funct. Anal.}, {\bf 97}, (1991).
\newblock 112-161.

\bibitem{Manujlov:94}
V.~M. Manuilov.
\newblock {Diagonalization of compact operators on Hilbert modules over
  W*-algebras of finite type (russ./engl.)}.
\newblock {\em Uspekhi Mat. Nauk}, {\bf 49}(2), (1994).
\newblock 159-160 / {\it Russ. Math. Surv.} {\bf 49}(1994), no.~2, 166-167.

\bibitem{Manujlov:94/4}
V.~M. Manuilov.
\newblock {On the eigenvalues of the perturbed Schr\"odinger operator with an
  irrational magnetic flow (russ./engl.)}.
\newblock {\em Funktional. Anal. i Prilozhen.}, {\bf 28}, (1994).
\newblock no.~2, 57-60 / {\it Funct. Anal. Appl.} {\bf 28}(1994), no.~2,
  120-122.

\bibitem{Manujlov:94/3}
V.~M. Manuilov.
\newblock {Representability of functionals and adjointability of operators on
  Hilbert C*-modules}.
\newblock preprint no.~1/94, Moscow State University, Dept. Mech. Math., Chair
  of Higher Geometry and Topology / submitted to Funktional. Anal. i Prilozh.,
  Moscow, 1994.

\bibitem{Manujlov:95/2}
V.~M. Manuilov.
\newblock {Diagonalization of compact operators on Hilbert modules over
  C*-algebras of real rank zero (engl.)}.
\newblock funct-an@babbage.sissa.it, preprint 9501008 / submitted to {\it
  Matem. Zametki}, Moskow, 1995.

\bibitem{Manujlov:94/2}
V.~M. Manuilov.
\newblock {Diagonalization of compact operators on Hilbert modules over finite
  W*-algebras (engl.)}.
\newblock {\em Annals Global. Anal. Geom.}, {\bf 13}, (1995).
\newblock 207-226.

\bibitem{Manujlov:95/1}
V.~M. Manuilov.
\newblock {Lusin's $C$-property is not valid for functional Hilbert modules}.
\newblock funct-an@babbage.sissa.it, preprint 9501004, 1995.

\bibitem{Manujlov:96}
V.~M. Manuilov.
\newblock {Diagonalizing operators over continuous fields of C*-algebras
  (engl.)}.
\newblock preprint, Moscow State Univ., Moscow, Russia, 1996.

\bibitem{Masani:59}
P.~R. Masani.
\newblock {Cram\'ers theorem on monotone matrix-valued functions and the Wold
  decomposition}.
\newblock In U.~Grenande, editor, {\em {Probability and Statistics - The Harald
  Cram\'er Volume}}. Almquist \& Wiksell, Stockholm, 1959.
\newblock 175-189.

\bibitem{Wiener/Masani:60}
P.~R. Masani.
\newblock {The prediction theory of multivariable stochastic processes, III.
  Unbounded spectral densities}.
\newblock {\em Acta Math.}, {\bf 104}, (1960).
\newblock 141-162.

\bibitem{Masani:62}
P.~R. Masani.
\newblock {Shift invariant spaces and prediction theory}.
\newblock {\em Acta Math.}, {\bf 107}, (1962).
\newblock 275-290.

\bibitem{Masani:66}
P.~R. Masani.
\newblock {Recent trends in multivariate prediction theory}.
\newblock In P.~R. Krishnaiah, editor, {\em {Multivariate Analysis, Proc. Int.
  Symp., Dayton, Ohio, June 1965}}. Academic Press, New York, 1966.
\newblock 351-381.

\bibitem{Mingo/Phillips:84}
J.~Mingo and W.~Phillips.
\newblock {Equivariant triviality theorems for Hilbert C*-modules}.
\newblock {\em Proc. Amer. Math. Soc.}, {\bf 91}, (1984).
\newblock 225-230.

\bibitem{Mingo:82/0}
J.~A. Mingo.
\newblock {\em {K-theory and multipliers of stable C*-algebras}}.
\newblock PhD thesis, Dalhousie Univ., Halifax, N.~S., Canada, 1982.

\bibitem{Mingo:82}
J.~A. Mingo.
\newblock {On the contractability of the unitary group of the Hilbert space
  over a C*-algebra}.
\newblock {\em Integral Equat. Operator Theory}, {\bf 5}, (1982).
\newblock 888-891.

\bibitem{Mingo:87}
J.~A. Mingo.
\newblock {K-theory and multipliers of stable C*-algebras}.
\newblock {\em Trans. Amer. Math. Soc.}, {\bf 299}, (1987).
\newblock 397-411.

\bibitem{Mingo:90}
J.~A. Mingo.
\newblock {Inner completely positive maps on von Neumann algebras}.
\newblock {\em Proc. Symp. Pure Math.}, {\bf 51}, (1990).
\newblock p. 2, 213-217.

\bibitem{Mishchenko:78}
A.~S. Mishchenko.
\newblock {The theory of elliptic operators over C*-algebras (russ./engl.)}.
\newblock {\em Dokl. Akad. Nauk SSSR}, {\bf 239}, (1978).
\newblock 1289-1291 / {\it Soviet Math. (Doklady)} {\bf 19}(1978), 512-515.

\bibitem{Mishchenko:79/1}
A.~S. Mishchenko.
\newblock {Banach algebras, pseudodifferential operators and their applications
  to K-theory (russ./engl.)}.
\newblock {\em Uspekhi Mat. Nauk}, {\bf 34}(6), (1979).
\newblock 67-79 / {\it Russ. Math. Surv.} {\bf 34}(1979), no. 6, 77-91.

\bibitem{Mishchenko:79/3}
A.~S. Mishchenko.
\newblock {C*-algebras and K-theory}.
\newblock {\em Lecture Notes Math.}, {\bf 763}, 1979.

\bibitem{Mishchenko:84}
A.~S. Mishchenko.
\newblock {Representations of compact groups on Hilbert modules over
  C*-algebras (russ./engl.)}.
\newblock {\em Trudy Mat. Inst. im. V. A. Steklova}, {\bf 166}, (1984).
\newblock 161-176 / {\it Proc. Steklov Inst. Math.} {\bf 166}(1986), 179-195.

\bibitem{Mishchenko:84/2}
A.~S. Mishchenko.
\newblock {\em {Vector bundles and their applications (russ.)}}.
\newblock Nauka, Moscow, USSR, 1984.

\bibitem{Mishchenko/Filippov:89}
A.~S. Mishchenko and O.~G. Filippov.
\newblock {The reduction of elliptical operators with almost periodical
  coefficients to operators on compact manifolds (russ./engl.)}.
\newblock {\em Vestn. Mosk. Univ., Ser. I: Mat.-Mekh.}, no. 5, 1989.
\newblock 78-81 / {\it Moscow Univ. Math. Bull.} {\bf 44}(1989), no. 5, 73-75.

\bibitem{Mishchenko/Fomenko:79/2}
A.~S. Mishchenko and A.~T. Fomenko.
\newblock {The index of elliptic operators over C*-algebras (russ./engl.)}.
\newblock {\em Izv. Akad. Nauk SSSR, Ser. Mat.}, {\bf 43}, (1979).
\newblock 831-859 / {\it Math. USSR - Izv.} {\bf 15}(1980), 87-112.

\bibitem{Mishchenko/Sharipov:83}
A.~S. Mishchenko and F.~Sharipov.
\newblock {Independence of spectra of elliptic operators with random
  coefficients (russ./engl.)}.
\newblock {\em Vestn. Mosk. Univ., Ser. I: Mat.-Mekh.}, no.6, (1983).
\newblock 51-56 / {\it Moscow Univ. Math. Bull.} {\bf 38}(1983), no. 6, 59-64.

\bibitem{Mishchenko/Solov'ov:77}
A.~S. Mishchenko and Yu.~P. Solov'ov.
\newblock {On infinite-dimensional representations of fundamental groups and on
  formulae of Hirzebruch type (russ./engl.)}.
\newblock {\em Dokl. Akad. Nauk SSSR}, {\bf 234}, (1977).
\newblock 761-764 / {\it Soviet Math. (Doklady)} {\bf 18}(1)(1977), 767-771.

\bibitem{Moore/Schochet:88}
C.~C. Moore and C.~Schochet.
\newblock {\em {Global analysis and foliated spaces}}.
\newblock Math. Sciences Research Institut Publication, no. 9.
  Springer--Verlag, New York, 1988.

\bibitem{Mori:88}
K.~Mori.
\newblock {Discrete series representations and K-theory of Hilbert C*-modules}.
\newblock {\em J. Fac. Sci. Technol., Kinki Univ.}, {\bf 24}, (1988).
\newblock 1-10.

\bibitem{Muhly/Renault/Williams:87}
P.~S. Muhly, J.~N. Renault, and D.~P. Williams.
\newblock {Equivalence and isomorphism for groupoid C*-algebras}.
\newblock {\em J. Oper. Theory}, {\bf 17}, (1987).
\newblock 3-22.

\bibitem{Muhly/Solel:95}
P.~S. Muhly and B.~Solel.
\newblock {Hilbert modules over operator algebras}.
\newblock {\em Memoirs Amer. Math. Soc.}, {\bf 559}, (1995).

\bibitem{Murphy:96}
G.~J. Murphy.
\newblock {Positive definite kernels and Hilbert C*-modules}.
\newblock preprint, University College Cork, Ireland / to appear in {\it Proc.
  Edinburgh Math. Soc.}, 1996.

\bibitem{Nagisa/Song:89}
M.~Nagisa and G.~Song.
\newblock {Inheritance of the solvability of the similarity problem within a
  C*-algebra and its C*-subalgebras}.
\newblock {\em Math. Jap.}, {\bf 34}, (1989).
\newblock 73-80.

\bibitem{Ng:95/2}
Chi-Keung Ng.
\newblock {Coactions and crossed products of Hopf C*-algebras, II: Hilbert
  C*-modules}.
\newblock preprint, Mathematical Institute, Oxford University, United Kingdom,
  1995.

\bibitem{Ng:95}
Chi-Keung Ng.
\newblock {Morita equivalences between fixed point algebras and crossed
  products}.
\newblock preprint, Mathematical Institute, Oxford University, United Kingdom,
  1995.

\bibitem{Ng:95/3}
Chi-Keung Ng.
\newblock {Morphisms of multiplicative unitaries}.
\newblock preprint, Mathematical Institute, Oxford University, United Kingdom,
  1995.

\bibitem{Ng:96}
Chi-Keung Ng.
\newblock {Discrete coactions on Hilbert C*-modules}.
\newblock {\em Math. Proc. Cambridge Philos. Soc.}, {\bf 119}, (1996).
\newblock 103-112.

\bibitem{Nishimura:90}
H.~Nishimura.
\newblock {Some connections between Boolean valued analysis and topological
  reduction theory for C*-algebras}.
\newblock {\em Z. Math. Logik Grundlagen Math.}, {\bf 36}, (1990).
\newblock 471-479.

\bibitem{Nishimura:93}
H.~Nishimura.
\newblock {A Boolean transfer principle from L*-algebras to AL*-algebras}.
\newblock {\em Math. Log. Q.}, {\bf 39}, (1993).
\newblock 241-250.

\bibitem{Olsen/Pedersen:89}
C.~L. Olsen and G.~K. Pedersen.
\newblock {Corona C*-algebras and their applications to lifting problems}.
\newblock {\em Math. Scand.}, {\bf 64}, (1989).
\newblock 63-68.

\bibitem{Ozawa:80}
M.~Ozawa.
\newblock {Hilbert B(H)-modules and stationary processes}.
\newblock {\em Kodai Math. J.}, {\bf 3}, (1980).
\newblock 26-39.

\bibitem{Ozawa:83}
M.~Ozawa.
\newblock {Boolean valued interpretation of Hilbert space theory}.
\newblock {\em J. Math. Soc. Japan}, {\bf 35}, (1983).
\newblock 609-627.

\bibitem{Ozawa:84}
M.~Ozawa.
\newblock {A classification of type I AW*-algebras and Boolean valued
  analysis}.
\newblock {\em J. Math. Soc. Japan}, {\bf 36}, (1984).
\newblock 589-608.

\bibitem{Ozawa:85/1}
M.~Ozawa.
\newblock {A transfer principle from von Neumann algebras to AW*-algebras}.
\newblock {\em J. London Math. Soc.}, {\bf 32}, (1985).
\newblock 141-148.

\bibitem{Ozawa:85/2}
M.~Ozawa.
\newblock {Nonuniqueness of the cardinality attached to homogeneous
  AW*-algebras}.
\newblock {\em Proc. Amer. Math. Soc.}, {\bf 93}, (1985).
\newblock 681-684.

\bibitem{Ozawa:86}
M.~Ozawa.
\newblock {Boolean valued analysis approach to the trace problem of
  AW*-algebras}.
\newblock {\em J. London Math. Soc.}, {\bf 33}, (1986).
\newblock 347-354.

\bibitem{Ozawa:90}
M.~Ozawa.
\newblock {Boolean valued interpretation of Banach space theory and module
  structures of von Neumann algebras}.
\newblock {\em Nagoya Math. J.}, {\bf 117}, (1990).
\newblock 1-36, (preprint no.12/85, Nagoya Univ., Japan).

\bibitem{Ozawa/Saito:86}
M.~Ozawa and K.~Sait\^o.
\newblock {Embedable AW*-algebras and regular completions}.
\newblock {\em J. London Math. Soc.}, {\bf 34}, (1986).
\newblock 511-523, (preprint no.6/85, Nagoya Univ. Japan).

\bibitem{Packer:86}
J.~A. Packer.
\newblock {K-theoretic invariants for C*-algebras associated to transformations
  and induced flows}.
\newblock {\em J. Funct. Anal.}, {\bf 67}, (1986).
\newblock 25-59.

\bibitem{Packer:87}
J.~A. Packer.
\newblock {C*-algebras generated by projective representations of the discrete
  Heisenberg group}.
\newblock {\em J. Oper. Theory}, (1987).
\newblock 41-66.

\bibitem{Packer:88/2}
J.~A. Packer.
\newblock {Flow equivalence for dynamical systems and the corresponding
  C*-algebras}.
\newblock {\em Oper. Theory: Adv. Appl.}, 28, 1988.
\newblock Birkh\"auser, Basel-Boston, Ma., 223-242.

\bibitem{Packer:88}
J.~A. Packer.
\newblock {Strong Morita equivalence for Heisenberg C*-algebras and the
  positive cones of their $K_o$-groups}.
\newblock {\em Canad. J. Math.}, {\bf 40}, (1988).
\newblock 833-864.

\bibitem{Packer/Raeburn:89}
J.~A. Packer and I.~Raeburn.
\newblock {Twisted crossed products of C*-algebras, I}.
\newblock {\em Math. Proc. Camb. Phil. Soc.}, {\bf 106}, (1989).
\newblock 293-311.

\bibitem{Packer/Raeburn:90}
J.~A. Packer and I.~Raeburn.
\newblock {Twisted crossed products of C*-algebras, II}.
\newblock {\em Math. Ann.}, {\bf 287}, (1990).
\newblock 595-612.

\bibitem{Papatriantafillou:94}
M.~H. Papatriantafillou.
\newblock {A reduction theorem for Hermitian structures on {\bf A}-bundles}.
\newblock {\em Bolletino d. Unione Matematica Italiana}, (7), {\bf 8 A},
  (1994).
\newblock 1-9.

\bibitem{Parker:88}
E.~M. Parker.
\newblock {The Brauer group of graded continuous trace C*-algebras}.
\newblock {\em Trans. Amer. Math. Soc.}, {\bf 308}, (1988).
\newblock 115-132.

\bibitem{Paschke:72}
W.~L. Paschke.
\newblock {\em {Hilbert B*-modules and completely bounded maps}}.
\newblock PhD thesis, University of Oregon, U.S.A., 1972.

\bibitem{Paschke:73}
W.~L. Paschke.
\newblock {Inner product modules over B*-algebras}.
\newblock {\em Trans. Amer. Math. Soc.}, {\bf 182}, (1973).
\newblock 443-468.

\bibitem{Paschke:74}
W.~L. Paschke.
\newblock {The double B-dual of an inner product module over a C*-algebra B}.
\newblock {\em Canad. J. Math.}, {\bf 26}, (1974).
\newblock 1272-1280.

\bibitem{Paschke:76}
W.~L. Paschke.
\newblock {Inner product modules arising from compact automorphism groups of a
  von Neumann algebra}.
\newblock {\em Trans. Amer. Math. Soc.}, {\bf 224}, (1976).
\newblock 87-102.

\bibitem{Paschke:77}
W.~L. Paschke.
\newblock {Integrable group actions on von Neumann algebras}.
\newblock {\em Math. Scand.}, {\bf 40}, (1977).
\newblock 234-248.

\bibitem{Paschke:85}
W.~L. Paschke.
\newblock {${\bf Z}_2$-equivariant K-theory}.
\newblock {\em Lecture Notes Math.}, {\bf 1132}, (1985).
\newblock Springer--Verlag, Berlin, pp. 362-373.

\bibitem{Pedersen:84}
G.~K. Pedersen.
\newblock {Multipliers in AW*-algebras}.
\newblock {\em Math. Z.}, {\bf 187}, (1984).
\newblock 23-24.

\bibitem{Pedersen:86}
G.~K. Pedersen.
\newblock {SAW*-algebras and corona C*-algebras. Contributions to
  non-commutative topology}.
\newblock {\em J. Operator Theory}, {\bf 15}, (1986).
\newblock 15-32.

\bibitem{Phillips/Raeburn:94}
J.~Phillips and I.~Raeburn.
\newblock {Twisted crossed products by coactions}.
\newblock {\em J. Austral. Math. Soc., Ser. A}, {\bf 56}, (1994).
\newblock 320-344.

\bibitem{Phillips:87}
N.~C. Phillips.
\newblock {Equivariant K-theory and freeness of group actions on C*-algebras}.
\newblock {\em Lecture Notes Math.}, {\bf 1274}, (1987).

\bibitem{Phillips:88/2}
N.~C. Phillips.
\newblock {Equivariant K-theory for proper actions and C*-algebras}.
\newblock {\em Contemp. Math.}, {\bf 70}, (1988).
\newblock 175-204.

\bibitem{Phillips:88/1}
N.~C. Phillips.
\newblock {Inverse limits of C*-algebras}.
\newblock {\em J. Oper. Theory}, {\bf 19}, (1988).
\newblock 159-195.

\bibitem{Phillips:89}
N.~C. Phillips.
\newblock {\em {Equivariant K-theory for proper actions}}.
\newblock Pitman Res. Notes Math., v. 178. Longman, London - New York, 1989.

\bibitem{Phillips:89/1}
N.~C. Phillips.
\newblock {Representable K-theory for $\sigma$-C*-algebras}.
\newblock {\em K-Theory}, {\bf 3}, (1989).
\newblock 441-478.

\bibitem{Pimsner:85}
M.~V. Pimsner.
\newblock {Range of traces on $K_o$ of reduced crossed products by free
  groups}.
\newblock {\em Lecture Notes Math.}, {\bf 1132}, (1985).

\bibitem{Pimsner:93}
M.~V. Pimsner.
\newblock {A class of C*-algebras generalizing both Cuntz-Krieger algebras and
  crossed products by {\bf Z}}.
\newblock preprint of the University of Pennsylvania, Dept. Math.,
  Philadelphia, PA, USA, 1993.

\bibitem{Pimsner/Popa:86}
M.~V. Pimsner and S.~Popa.
\newblock {Entropy and index for subfactors}.
\newblock {\em Ann. scient. \'Ec. Norm. Sup. $4^e$ s\'erie}, {\bf 19}, (1986).
\newblock 57-106.

\bibitem{Pimsner/Voiculescu:80}
M.~V. Pimsner and D.~Voiculescu.
\newblock {Exact sequences for K-groups and EXT-groups of certain crossed
  product C*-algebras}.
\newblock {\em J. Oper. Theory}, {\bf 4}, (1980).
\newblock 201-210.

\bibitem{Pincket:86}
J.~Pincket.
\newblock {Hilbert modules over an arbitrary C*-algebra}.
\newblock {\em Bull. Soc. Math. Belg., S\'er. B}, {\bf 38}, (1986).
\newblock 176-186.

\bibitem{Plymen:86}
R.~J. Plymen.
\newblock {Strong Morita equivalence, spinors and symplectic spinors}.
\newblock {\em J. Oper. Theory}, {\bf 16}, (1986).
\newblock 305-324.

\bibitem{Plymen:87}
R.~J. Plymen.
\newblock {The reduced C*-algebra of the p-adic group GL(n)}.
\newblock {\em J. Func. Anal.}, {\bf 72}, (1987).
\newblock 1-12.

\bibitem{Plymen:90}
R.~J. Plymen.
\newblock {Equivalence bimodules in the representation theory of reductive
  groups}.
\newblock {\em Proc. Symp. Pure Math.}, {\bf 51-1}, (1990).
\newblock 267-272.

\bibitem{Putnam:85}
I.~F. Putnam.
\newblock {Strong Morita equivalence of Denjoy C*-algebras}.
\newblock {\em C. R. Math. Rep. Acad. Sci. Canada}, {\bf 7}, (1985).
\newblock 121-125.

\bibitem{Putnam:88}
I.~F. Putnam.
\newblock {Strong Morita equivalence for the Denjoy C*-algebras}.
\newblock {\em Can. Math. Bull.}, {\bf 31}, (1988).
\newblock 439-447.

\bibitem{Putnam/Schmidt/Skau:86}
I.~F. Putnam, K.~Schmidt, and C.~Skau.
\newblock {C*-algebras associated with Denjoy homeomorphisms of the circle}.
\newblock {\em J. Oper. Theory}, {\bf 16}, (1986).
\newblock 99-126.

\bibitem{Qing:95}
Lin Qing.
\newblock {Cut-down method in the inductive limit decomposition of
  non-commutative tori, III: A complete answer in 3-dimension}.
\newblock preprint, 1995.

\bibitem{Quigg:95}
J.~C. Quigg.
\newblock {Full and reduced C*-coactions}.
\newblock {\em Math. Proc. Camb. Phil. Soc.}, {\bf 116}, (1995).
\newblock 435-450.

\bibitem{Quigg/Spielberg}
J.~C. Quigg and J.~Spielberg.
\newblock {Regularity and hyporegularity in C*-dynamical systems}.
\newblock {\em Houston Math. J.}, {\bf 18}, (1992).
\newblock 139-152.

\bibitem{Raeburn:81}
I.~Raeburn.
\newblock {On the Picard group of a continuous trace C*-algebra}.
\newblock {\em Trans. Amer. Math. Soc.}, {\bf 263}, (1981).
\newblock 183-205.

\bibitem{Raeburn:88}
I.~Raeburn.
\newblock {Induced C*-algebras and a symmetric imprimitivity theorem}.
\newblock {\em Math. Ann.}, {\bf 280}, (1988).
\newblock 369-387.

\bibitem{Raeburn:95}
I.~Raeburn.
\newblock {Induced C* -algebras and a symmetric imprimitivity theorem}.
\newblock preprint, Sch. Math., Univ. NSW, Kensington, NSW 2033, Aust., to
  appear in Math. Ann., 1995.

\bibitem{Raeburn/Williams:85}
I.~Raeburn and D.~P. Williams.
\newblock {Pull-backs of C*-algebras and crossed products by certain diagonal
  actions}.
\newblock {\em Trans. Amer. Math. Soc.}, {\bf 287}, (1985).
\newblock 755-777.

\bibitem{Raeburn/Williams:93}
I.~Raeburn and D.~P. Williams.
\newblock {Dixmier-Douady classes of dynamical systems and crossed products}.
\newblock {\em Can. J. Math.}, {\bf 45}, (1993).
\newblock 1032-1066.

\bibitem{Rieffel:74/2}
M.~A. Rieffel.
\newblock {Induced representations of C*-algebras}.
\newblock {\em Adv. Math.}, {\bf 13}, (1974).
\newblock 176-257.

\bibitem{Rieffel:74/1}
M.~A. Rieffel.
\newblock {Morita equivalence for C*-algebras and W*-algebras}.
\newblock {\em J. Pure Applied Alg.}, {\bf 5}, (1974).
\newblock 51-96.

\bibitem{Rieffel:76}
M.~A. Rieffel.
\newblock {Strong Morita equivalence of certain transformation group
  C*-algebras}.
\newblock {\em Math. Ann.}, {\bf 222}, (1976).
\newblock 7-22.

\bibitem{Rieffel:81}
M.~A. Rieffel.
\newblock {C*-algebras associated with irrational rotations}.
\newblock {\em Pac. J. Math.}, {\bf 93}, (1981).
\newblock 415-429.

\bibitem{Rieffel:82/2}
M.~A. Rieffel.
\newblock {Applications of strong Morita equivalence to transformation group
  C*-algebras}.
\newblock {\em Proc. Symp. Pure Math. Amer. Math. Soc.}, {\bf 38}(1), (1982).
\newblock 299-310.

\bibitem{Rieffel:82/1}
M.~A. Rieffel.
\newblock {Morita equivalence for operator algebras}.
\newblock {\em Proc. Symp. Pure Math. Amer. Math. Soc.}, {\bf 38}(1), (1982).
\newblock 285-298.

\bibitem{Rieffel:83/1}
M.~A. Rieffel.
\newblock {Dimension and stable rank in the K-theory of C*-algebras}.
\newblock {\em Proc. London Math. Soc.}, {\bf 47}, (1983).
\newblock 285-302.

\bibitem{Rieffel:83/2}
M.~A. Rieffel.
\newblock {The cancellation theorem for projective modules over irrational
  rotation algebras}.
\newblock {\em Proc. London Math. Soc. (3)}, {\bf 47}, (1983).
\newblock 285-302.

\bibitem{Rieffel:85}
M.~A. Rieffel.
\newblock {''Vector bundles'' over higher dimensional non-commutative tori}.
\newblock {\em Lecture Notes Math.}, {\bf 1132}, (1985).
\newblock Springer--Verlag, Berlin, pp. 456-467.

\bibitem{Rieffel:87/1}
M.~A. Rieffel.
\newblock {Non-stable K-theory and non-commutative tori}.
\newblock {\em Contemp. Math.}, {\bf 62}, (1987).
\newblock 267-279.

\bibitem{Rieffel:87/2}
M.~A. Rieffel.
\newblock {The homotopy groups of the unitary groups of non-commutative tori}.
\newblock {\em J. Oper. Theory}, {\bf 17}, (1987).
\newblock 237-254.

\bibitem{Rieffel:88/1}
M.~A. Rieffel.
\newblock {Projective modules over higher-dimensional non-commutative tori}.
\newblock {\em Canad. J. Math.}, {\bf 40}, (1988).
\newblock 257-338.

\bibitem{Rieffel:90/2}
M.~A. Rieffel.
\newblock {Critical points of Yang--Mills for non-commutative tori}.
\newblock {\em J. Differential Geom.}, {\bf 31}, (1990).
\newblock 535-546.

\bibitem{Rieffel:90/1}
M.~A. Rieffel.
\newblock {Non-commutative tori -- a case study of non-commutative
  differentiable manifolds}.
\newblock {\em Contemp. Math.}, {\bf 105}, (1990).

\bibitem{Rieffel:91}
M.~A. Rieffel.
\newblock {\em {Proper actions of groups on C*-algebras}}, pages 141--182.
\newblock PM v.{\bf 84}. Mappings of Operator Algebras, ed.: H. Araki, R. V.
  Kadison, Birkh\"auser--Verlag, Boston- Basel - Berlin, 1991.

\bibitem{Rieffel:93}
M.~A. Rieffel.
\newblock {Deformation quantization for actions of ${\bf R}^d$}.
\newblock {\em Memoirs Amer. Math. Soc.}, {\bf 506}, (1993).

\bibitem{Rieffel:93/2}
M.~A. Rieffel.
\newblock {K-groups of C*-algebras deformed by actions of ${\bf R}^d$}.
\newblock {\em J. Funct. Anal.}, {\bf 116}, (1993).
\newblock 199-214.

\bibitem{Rosenberg:83}
J.~Rosenberg.
\newblock {C*-algebras, positive scalar curvature, and the Novikov conjecture}.
\newblock {\em Publ. Math. I.H.E.S.}, {\bf 58}, (1983).
\newblock 197-212.

\bibitem{Rosenberg:88}
J.~Rosenberg.
\newblock {K-theory of group C*-algebras, foliation C*-algebras and crossed
  products}.
\newblock {\em Contemp. Math.}, {\bf 70}, (1988).
\newblock pp. ?

\bibitem{Rosenberg:90}
J.~Rosenberg.
\newblock {K and KK: Topology and operator algebras}.
\newblock {\em Proc. Symp. Pure Math.}, {\bf 51-1}, (1990).
\newblock 445-480.

\bibitem{Rosenberg:64}
M.~Rosenberg.
\newblock {The square-integrability of matrix-valued functions with respect to
  a non-negative hermitian measure}.
\newblock {\em Duke Math. J.}, {\bf 31}, (1964).
\newblock 291-298.

\bibitem{Saito:71}
K.~Sait{\^o}.
\newblock {On the embedding as a double commutant in a type I AW*-algebra, I}.
\newblock {\em T{\^o}hoku Math. J.}, {\bf 23}, (1971).
\newblock 541-557.

\bibitem{Saito:74}
K.~Sait{\^o}.
\newblock {On the embedding as a double commutant in a type I AW*-algebra, II}.
\newblock {\em T{\^o}hoku Math. J.}, {\bf 26}, (1974).
\newblock 333-339.

\bibitem{Saito:78}
K.~Sait{\^o}.
\newblock {AW*-algebras with monotone convergence property and type III,
  non-W*, AW*-factors}.
\newblock {\em Lect. Notes Math.}, {\bf 650}, (1978).
\newblock 131-134.

\bibitem{Saito:79}
K.~Saito.
\newblock {AW*-algebras with monotone convergence property and examples by
  Takenouchi and Dyer}.
\newblock {\em T{\^o}hoku Math. J.}, {\bf 31}, (1979).
\newblock 31-40.

\bibitem{Saito:94}
K.~Sait{\^o}.
\newblock {Wild, type III, monotone complete, simple C*-algebras indexed by
  cardinal numbers}.
\newblock {\em J. London Math. Soc.(2)}, {\bf 49}, (1994).
\newblock 543-554.

\bibitem{Salehi:65}
H.~Salehi.
\newblock {\em {The prediction theory of multivariate stochastic processes with
  continuous time}}.
\newblock PhD thesis, Indiana Univ., Bloomington, U.S.A., 1965.

\bibitem{Salehi:66}
H.~Salehi.
\newblock {A factorization algorithm for q$\times$q matrix-valued functions on
  the real line R}.
\newblock {\em Trans. Amer. Math. Soc.}, {\bf 124}, (1966).
\newblock 468-470.

\bibitem{Salehi:67}
H.~Salehi.
\newblock {On the growth of a q-variate stationary stochastic process}.
\newblock {\em Z. Wahrscheinlichkeitstheorie und verw. Gebiete}, {\bf 8},
  (1967).
\newblock 140-147.

\bibitem{Sauvageot:89}
J.-L. Sauvageot.
\newblock {Tangent bimodule and locality for dissipative operators on
  C*-algebras}.
\newblock {\em Lecture Notes Math.}, {\bf 1396}, (1989).
\newblock Springer--Verlag, Berlin, pp. ?

\bibitem{Saworotnow:68}
P.~P. Saworotnow.
\newblock {A generalized Hilbert space}.
\newblock {\em Duke Math. J.}, {\bf 35}, (1968).
\newblock 191-197.

\bibitem{Saworotnow:83}
P.~P. Saworotnow.
\newblock {Generalized positive definite functions and stationary processes}.
\newblock In V.~Mandrekar and H.~Salehi, editors, {\em {Prediction Theory and
  Harmonic Analysis, The Pesi Masani volume}}. North-Holland Publishing Comp.,
  Amsterdam, 1983.
\newblock 329-344.

\bibitem{Scedrow/Scowcroft:88}
A.~Scedrow and P.~Scowcroft.
\newblock {Decomposition of finitely generated modules over C(X): sheaf
  semantics and a decision procedure}.
\newblock {\em Math. Proc. Cambridge Philos. Soc.}, {\bf 103}, (1988).
\newblock 257-268.

\bibitem{Schochet:94}
C.~Schochet.
\newblock {Equivariant KK-theory for inverse limits of $G$-C*-algebras}.
\newblock {\em J. Austral. Math. Soc. (Series A)}, {\bf 56}, (1994).
\newblock 183-211.

\bibitem{Schroeder:93}
H.~Schr\"oder.
\newblock {\em {K-theory for real C*-algebras and applications}}, volume~{\bf
  290} of {\em Pitman Res. Notes in Math. Sci.}
\newblock Longman Scientific \& Technical, Harlow, England, 1993.

\bibitem{Sharipov:85/1}
F.~Sharipov.
\newblock {Independence of the spectrum of an elliptic operator over a
  C*-algebra (russ./engl.)}.
\newblock {\em Vestn. Mosk. Univ., Ser. I: Mat.-Mekh.}, no. 1, (1985).
\newblock 87-89 / {\it Moscow Univ. Math. Bull.} {\bf 40}(1985), no. 1, 96-99.

\bibitem{Sharipov:85/2}
F.~Sharipov.
\newblock {Representation of the C*-algebra $End_A^*(l_2(A))$ (russ.)}.
\newblock In {\em {Problems in mathematical analysis and its applications}}.
  Gos. Univ. Samarkand, USSR, (1985).
\newblock pp. 89-93.

\bibitem{Sharipov/Zhuraev:86}
F.~Sharipov and Yu.~I. Zhuraev.
\newblock {On the index of a Fredholm operator in a Hilbert C*-module (russ.)}.
\newblock In {\em {Problems in math. anal. and its appl.}} Samarkand. Gos.
  Univ., Samarkand, (1986).
\newblock pp. 37-41.

\bibitem{Sharipov/Zhuraev:88}
F.~Sharipov and Yu.~I. Zhuraev.
\newblock {Index of elliptic operators over a C*-algebra (russ.)}.
\newblock In {\em {Problems in multidimensional differential geometry and its
  applications}}. Samarkand. Gos. Univ., Samarkand, (1988).
\newblock pp. 47-52.

\bibitem{Shen:82}
Nien-Tsu Shen.
\newblock {\em {Embeddings of Hilbert bimodules}}.
\newblock PhD thesis, Purdue University, West Lafayette, USA, 1982.

\bibitem{Sheu:87/2}
A.~J.~L. Sheu.
\newblock {A cancellation theorem for modules over the group C*-algebras of
  certain nilpotent Lie groups}.
\newblock {\em Canad. J. Math.}, {\bf 39}, (1987).
\newblock 365-427.

\bibitem{Sheu:87/1}
A.~J.~L. Sheu.
\newblock {Classification for projective modules over the unitized group
  C*-algebras for certain solvable Lie-groups}.
\newblock {\em J. Oper. Theory}, {\bf 18}, (1987).
\newblock 33-40.

\bibitem{Skandalis:84}
G.~Skandalis.
\newblock {Some remarks on Kasparov's theory}.
\newblock {\em J. Func. Anal.}, {\bf 56}, (1984).
\newblock 337-347.

\bibitem{Skandalis:88}
G.~Skandalis.
\newblock {Une notion de nucle\'arit\'e en K-th\'eorie (d'apr\`es J. Cuntz)}.
\newblock {\em K-theory}, {\bf 1}, (1988).
\newblock 549-573.

\bibitem{Skandalis:91}
G.~Skandalis.
\newblock {Kasparov's bivariant K-theory and applications}.
\newblock {\em Expo. Math.}, {\bf 9}, (1991).
\newblock 193-250.

\bibitem{Sunouchi:71}
C.~Sunouchi.
\newblock {A generalization of Schatten--von Neumann--Dixmier theorem for type
  I AW*-algebras}.
\newblock {\em T{\^o}hoku Math. J.}, {\bf 23}, (1971).
\newblock 727-734.

\bibitem{Swan:62}
R.~G. Swan.
\newblock {Vector bundles and projective modules}.
\newblock {\em Trans. Amer. Math. Soc.}, {\bf 105}, (1962).
\newblock 264-277.

\bibitem{Swan:77}
R.~G. Swan.
\newblock {Topological examples of projective modules}.
\newblock {\em Trans. Amer. Math. Soc.}, {\bf 230}, (1977).
\newblock 201-234.

\bibitem{Swan:87}
R.~G. Swan.
\newblock {Vector bundles, projective modules and the K-theory of spheres}.
\newblock In {\em {Algebraic topology and algebraic K-theory}}. Princeton Univ.
  Press, Princeton, NJ, 1987.
\newblock (Ann. of Math. Stud. v.113), pp. 432-522.

\bibitem{Takahashi:79/2}
A.~Takahashi.
\newblock {A duality between Hilbert modules and fields of Hilbert spaces}.
\newblock {\em Rev. Colomb. Mat.}, {\bf 13}, (1979).
\newblock 93-120.

\bibitem{Takahashi:79/1}
A.~Takahashi.
\newblock {Hilbert modules and their representations}.
\newblock {\em Rev. Colomb. Mat.}, {\bf 13}, (1979).
\newblock 1-38.

\bibitem{Takahashi:71}
A.~O. Takahashi.
\newblock {\em {Fields of Hilbert modules}}.
\newblock PhD thesis, Tulane Univ., New Orleans, USA, 1971.

\bibitem{Takemoto:73/1}
H.~Takemoto.
\newblock {On a characterization of AW*-modules and a representation of Gelfand
  type of non-commutative operator algebras}.
\newblock {\em Michigan Math. J.}, {\bf 20}, (1973).
\newblock 115-127.

\bibitem{Takemoto:75}
H.~Takemoto.
\newblock {Decomposable operators in continuous fields of Hilbert spaces}.
\newblock {\em T{\^o}hoku Math. J.}, {\bf 27}, (1975).
\newblock 413-435.

\bibitem{Takemoto:76}
H.~Takemoto.
\newblock {On the weakly continuous constant field of Hilbert space and its
  application to the reduction theory of von Neumann algebras}.
\newblock {\em T{\^o}hoku Math. Journal}, {\bf 28}, (1976).
\newblock 479-496.

\bibitem{Takenouchi:78}
O.~Takenouchi.
\newblock {A non-W*, AW*-factor}.
\newblock {\em Lect. Notes Math.}, {\bf 650}, (1978).
\newblock 135-139.

\bibitem{Takesaki:60}
M.~Takesaki.
\newblock {On the Hahn--Banach type theorem and the Jordan decomposition of
  module linear mapping over some operator algebras}.
\newblock {\em Kodai Math. Sem. Reports}, {\bf 12}, (1960).
\newblock 1-10.

\bibitem{Takeuti:83/2}
G.~Takeuti.
\newblock {C*-algebras and Boolean valued analysis}.
\newblock {\em Japan. J. Math.}, {\bf 9}, (1983).
\newblock 207-246.

\bibitem{Takeuti:83/1}
G.~Takeuti.
\newblock {Von Neumann algebras and Boolean valued analysis}.
\newblock {\em J. Math. Soc. Japan}, {\bf 35}, (1983).
\newblock 1-21.

\bibitem{Thomsen:88}
K.~Thomsen.
\newblock {\em {Hilbert C*-modules, KK-theory and C*-extensions}}.
\newblock Various Publications Series v.38. Aarhus Universitet, Matematisk
  Institut, Aarhus, 1988.

\bibitem{Tomiyama:87}
J.~Tomiyama.
\newblock {\em {Invitation to C*-algebras and topological dynamics}}.
\newblock World Scientific, Singapore, 1987.

\bibitem{Trofimov:86}
V.~A. Trofimov.
\newblock {Reflexivity of Hilbert modules over the algebra of compact operators
  with adjoint identity (russ./engl.)}.
\newblock {\em Vestn. Mosk. Univ., Ser. I: Mat.-Mekh.}, no.5, (1986).
\newblock 60-64 / {\it Moscow Univ. Math. Bull.} {\bf 41}(1986), no.5, 51-55.

\bibitem{Trofimov:87/2}
V.~A. Trofimov.
\newblock {Reflexive and self-dual Hilbert modules over some C*-algebras
  (russ./engl.)}.
\newblock {\em Uspekhi Mat. Nauk}, {\bf 42}, (1987).
\newblock 247-248 / {\it Russian Math. Surveys} {\bf 42}(1987), 303-304.

\bibitem{Trofimov:87/1}
V.~A. Trofimov.
\newblock {Reflexivity of Hilbert modules over splittable extensions of the
  algebra of compact operators (russ.)}.
\newblock In {\em {Geometry and the theory of singularities in nonlinear
  equations}}, Voronesh (USSR), 1987. Voronesh. Gos. Univ.
\newblock pp.164-170.

\bibitem{Trofimov:87/3}
V.~A. Trofimov.
\newblock {Self-dual and reflexive Hilbert modules over C*-algebras (russ.)}.
\newblock In {\em {Bakinskaya mezhdunarodnaya topologi{\^c}eskaya konferenciya,
  proceedings, part 2}}, Baku (USSR), 1987.
\newblock p. 296.

\bibitem{Trofimov:87/4}
V.~A. Trofimov.
\newblock {\em {The structure of Hilbert modules over topological spaces and
  over operator algebras (russ.)}}.
\newblock PhD thesis, Moscow State University ''M. V. Lomonosov'', Moscow
  (USSR), (1987).

\bibitem{Troitskij:85/3}
E.~V. Tro\H{\i}tsky.
\newblock {A connection between complex and operator topological equivariant
  K-theories (russ./engl.)}.
\newblock {\em Uspekhi Mat. Nauk}, {\bf 40}, (1985).
\newblock 227-228 / {\it Russian Math. Surv.} {\bf 40}(1985), no. 4, 243-244.

\bibitem{Troitskij:85/4}
E.~V. Tro\H{\i}tsky.
\newblock {A theorem on the index of equivariant C*-elliptical operators
  (russ./engl.)}.
\newblock {\em Dokl. Akad. Nauk USSR}, {\bf 282}, (1985).
\newblock 1059-1061 / {\it Russian Math. Surv.} {\bf 31}(1985), 558-560.

\bibitem{Troitskij:85/1}
E.~V. Tro\H{\i}tsky.
\newblock {The index theorem for equivariant C*-elliptical operators (russ.)}.
\newblock {\em Dokl. Akad. Nauk USSR}, {\bf 282}, (1985).
\newblock 1059-1061.

\bibitem{Troitskij:85/2}
E.~V. Tro\H{\i}tsky.
\newblock {The representation space of the K-functor related to a C*- algebra
  (russ./engl.)}.
\newblock {\em Vestn. Mosk. Univ., Ser. I: Mat.-Mekh.}, no. 1, (1985).
\newblock 96-98 / {\it Moscow Univ. Math. Bull.} {\bf 40}(1985), no. 1,
  111-115.

\bibitem{Troitskij:86/2}
E.~V. Tro\H{\i}tsky.
\newblock {An equivariant index theorem with C*-elliptic operators
  (russ./engl.)}.
\newblock {\em Izv. Akad. Nauk SSSR, Ser. Mat.}, {\bf 50}, (1986).
\newblock 849-865 / {\it Math. USSR - Izv.} {\bf 29}(1986), 207-224.

\bibitem{Troitskij:86}
E.~V. Tro\H{\i}tsky.
\newblock {Contractability of the full general linear group of the Hilbert
  C*-module $l_2(A)$ (russ./engl.)}.
\newblock {\em Funktsional. Anal. i Prilozh.}, {\bf 20}(4), (1986).
\newblock 58-64 / {\it Funct. Anal. Appl.} {\bf 20}(1986), 301-307.

\bibitem{Troitskij:86/3}
E.~V. Tro\H{\i}tsky.
\newblock {Homotopic triviality of the general linear group of a Hilbert module
  (russ.)}.
\newblock In V.~V. Kozlov and A.~T. Fomenko, editors, {\em Geom., diff.
  equations and mechanics}, Moscow State Univ., Mekh.-Mat. Fakulty, Moscow,
  1986.
\newblock pp. 128-134.

\bibitem{Troitskij:87}
E.~V. Tro\H{\i}tsky.
\newblock {The index of equivariant elliptic operators over C*-algebras}.
\newblock {\em Annals Global Anal. Geom.}, {\bf 5}, (1987).
\newblock 3-22.

\bibitem{Troitskij:88}
E.~V. Tro\H{\i}tsky.
\newblock {An exact K-cohomology C*-index formula, I: Thom isomorphism and
  topological index (russ./engl.)}.
\newblock {\em Vestn. Mosk. Univ., Ser. I: Mat.-Mekh.}, no. 2, (1988).
\newblock 83-85 / {\it Moscow Univ. Math. Bull.} {\bf 43}(1988), no. 2, 57-60.

\bibitem{Troitskij:89}
E.~V. Tro\H{\i}tsky.
\newblock {An exact K-cohomology C*-index formula, II: An index theorem and its
  applications (russ./engl.)}.
\newblock {\em Uspekhi Mat. Nauk}, {\bf 44}, (1989).
\newblock 213-214 / {\it Russian Math. Surv.} {\bf 44}(1989), 259-261.

\bibitem{Troitskij:91}
E.~V. Tro\H{\i}tsky.
\newblock {Lefschetz numbers of C*-complexes}.
\newblock {\em Lect. Notes Math.}, {\bf 1474}, (1991).
\newblock 193-206.

\bibitem{Troitskij:92}
E.~V. Tro\H{\i}tsky.
\newblock {An exact formula for the index of an equivariant C*-elliptic
  operator (russ./engl.)}.
\newblock {\em Tr. Mat. Inst. Steklova}, {\bf 193}, (1992).
\newblock 178-182 / {\it Proc. Steklov Inst. Math. {\bf 193}(1993), 197-201}.

\bibitem{Troitskij:93/2}
E.~V. Tro\H{\i}tsky.
\newblock {Traces, C*-elliptic complexes and higher even cyclic homology}.
\newblock {\em Vestn. Mosk. Univ., Ser. I: Mat.-Mekh.}, no. 5, (1993).
\newblock 36-39 / {\it Moscow Univ. Math. Bull.} {\bf 48}(1993), no. 5.

\bibitem{Troitskij:94/2}
E.~V. Tro\H{\i}tsky.
\newblock {An averaging theorem in Hilbert C*-modules and operators possessing
  an adjoint (russ./engl.)}.
\newblock {\em Funktsional. Anal. i Prilozh.}, {\bf 28}(3), (1994).
\newblock 88-92 / {\it Funct. Anal. Appl.} {\bf 28}(1994), no.~3, 220-223.

\bibitem{Troitskij:94}
E.~V. Tro\H{\i}tsky.
\newblock {Some aspects of geometry of operators in Hilbert modules}.
\newblock preprint, Ruhr-Universit\"at Bochum, Fakult\"at f\"ur Mathematik,
  Bericht-Nr. 173, January, 1994.

\bibitem{Troitskij:95}
E.~V. Tro\H{\i}tsky.
\newblock {Operators without adjoint and Kuiper's theorem for Hilbert modules}.
\newblock {\rm Universit\"at Heidelberg, Mathematisches Institut, preprint
  series of Forschergruppe ''Topologie und Nichtkommutative Geometrie'', no.
  113, Februar}, 1995.

\bibitem{Troitskij:96}
E.~V. Tro\H{\i}tsky.
\newblock {Kuiper's theorem for Hilbert modules: the general case}.
\newblock Max-Planck-Institut f\"ur Math., Bonn, preprint MPI 96-16, 1996.

\bibitem{Truong-Van:81}
B.~Truong-Van.
\newblock {Une g\'en\'eralisation du th\'eor\`eme de Kolmogorov-Aronszajn.
  Processus V-born\'es q-dimensionelles: domaine spectral}.
\newblock {\em Ann. Inst. Henri Poincar\'e}, {\bf 17}, (1981).
\newblock 31-49.

\bibitem{Tsui:96}
Sze-Kai Tsui.
\newblock {Completely positive module maps and completely positive extreme
  maps}.
\newblock {\em Proc. Amer. Math. Soc.}, {\bf 124}, (1996).
\newblock 437-445.

\bibitem{Umegaki:55}
H.~Umegaki.
\newblock {Positive definite function and a direct product Hilbert space}.
\newblock {\em T{\^o}hoku Math. J.}, {\bf 7}, (1955).
\newblock 206-211.

\bibitem{Varela:74}
J.~Varela.
\newblock {Sectional representations of Banach modules}.
\newblock {\em Math. Z.}, {\bf 139}, (1974).
\newblock 55-61.

\bibitem{Vincent--Smith:67}
G.~Vincent-Smith.
\newblock {The Hahn--Banach theorem for modules}.
\newblock {\em Proc. London Math. Soc.}, {\bf 17}, (1967).
\newblock 72-90.

\bibitem{Walters:94}
S.~G. Walters.
\newblock {Strong Morita equivalence for the quasi-rotation C*-algebra}.
\newblock {\em J. Operator Theory}, {\bf 31}, (1994).
\newblock 327-349.

\bibitem{Walters:95}
S.~G. Walters.
\newblock {Projective modules over the non-commutative sphere}.
\newblock {\em J. London Math. Soc. (2)}, {\bf 51}, (1995).
\newblock 589-602.

\bibitem{Wang:89}
X.~Wang.
\newblock {\em {The C*-algebras of a class of solvable Lie groups}}.
\newblock Pitman Res. Notes Math. v.199. Longman Scientific, Harlow, 1989.

\bibitem{Watatani:90}
Y.~Watatani.
\newblock {Index for C*-subalgebras}.
\newblock {\em Memoirs Amer. Math. Soc.}, {\bf 424}, (1990).

\bibitem{Weaver:95/2}
N.~Weaver.
\newblock {Deformations of von Neumann algebras}.
\newblock preprint, Univ. of California at Santa Barbara, U.S.A., 1995.

\bibitem{Weaver:95/1}
N.~Weaver.
\newblock {Full C*-dynamical systems}.
\newblock preprint, Univ. of California at Santa Barbara, U.S.A., 1995.

\bibitem{Wegge--Olsen:89}
N.~E. Wegge-Olsen.
\newblock {Introduction to operator algebra K-theory and generalized index
  theory}.
\newblock preprints no.1a+1b, K{\o}benhavns Univ., Matematisk Institut, 1989.

\bibitem{Wegge--Olsen:93}
N.~E. Wegge-Olsen.
\newblock {\em {K-theory and C*-algebras -- a friendly approach}}.
\newblock Oxford University Press, Oxford, 1993.

\bibitem{Weidner:87}
J.~Weidner.
\newblock {\em {Topological invariants for generalized operator algebras}}.
\newblock PhD thesis, Univ. Heidelberg, Heidelberg, FRG, 1987.

\bibitem{Weidner:89/1}
J.~Weidner.
\newblock {KK-groups for generalized operator algebras,I}.
\newblock {\em K-Theory}, {\bf 3}, (1989).
\newblock 57-77.

\bibitem{Weidner:89/2}
J.~Weidner.
\newblock {KK-groups for generalized operator algebras,II}.
\newblock {\em K-Theory}, {\bf 3}, (1989).
\newblock 79-98.

\bibitem{Wickstead:82}
A.~W. Wickstead.
\newblock {Stone-algebra-valued measures: Integration of vector-valued
  functions and Radon-Nikodym type theorems}.
\newblock {\em Proc. Lond. Math. Soc., III. Ser.}, {\bf 45}, (1982).
\newblock 193-226.

\bibitem{Widom:56}
H.~Widom.
\newblock {Embedding in algebras of type I}.
\newblock {\em Duke Math. J.}, {\bf 23}, (1956).
\newblock 309-324.

\bibitem{Wiener/Masani:57}
N.~Wiener and P.~Masani.
\newblock {The prediction theory of multivariable stochastic processes, I. The
  regularity condition}.
\newblock {\em Acta Math.}, {\bf 98}, (1957).
\newblock 111-150.

\bibitem{Wiener/Masani:58}
N.~Wiener and P.~Masani.
\newblock {The prediction theory of multivariable stochastic processes, II. The
  linear predictor}.
\newblock {\em Acta Math.}, {\bf 99}, (1958).
\newblock 93-137.

\bibitem{Wittstock:81}
G.~Wittstock.
\newblock {Ein operatorwertiger Hahn--Banach--Satz}.
\newblock {\em J. Funct. Anal.}, {\bf 40}, (1981).
\newblock 127-150.

\bibitem{Wittstock:84}
G.~Wittstock.
\newblock {Extensions of completely bounded C*-module homomorphisms}.
\newblock In {\em Operator algebras and group representations, II, (Monographs
  and Studies in Math. v.18), Proc. Int. Conf. Neptun(Rom.), 1980}. Pitman Adv.
  Publishing Progr., Boston, 1984.
\newblock pp. ?

\bibitem{Woronowicz:91}
S.~L. Woronowicz.
\newblock {Unbounded elements affiliated with C*-algebras and non-compact
  quantum groups}.
\newblock {\em Commun. Math. Phys.}, {\bf 136}, (1991).
\newblock 399-432.

\bibitem{Woronowicz/Napiorkowski:92}
S.~L. Woronowicz and K.~Napi\'orkowski.
\newblock {Operator theory in the C*-algebra framework}.
\newblock {\em Rep. Math. Phys.}, (1992).
\newblock 353-371.

\bibitem{Wright:69/2}
J.~D.~M. Wright.
\newblock {A spectral theorem for normal operators an a Kaplansky--Hilbert
  module}.
\newblock {\em Proc. London Math. Soc.}, {\bf 19}, (1969).
\newblock 258-268.

\bibitem{Xu:91}
P.~Xu.
\newblock {Morita equivalent symplectic groupoids}.
\newblock In {\em {Symplectic geometry, groupoids and integrable systems
  (Berkeley, CA, 1989)}}, New York, 1991. Springer-Verlag.
\newblock (Math. Sci. Res. Inst. Publ. v. 20).

\bibitem{Xu:92}
P.~Xu.
\newblock {Morita equivalence and symplectic realizations of Poisson
  manifolds}.
\newblock {\em Ann. Sci. Ec. Norm. Super., Iv.~Ser.}, {\bf 25}, (1992).
\newblock 307-333.

\bibitem{Yamagami:84}
S.~Yamagami.
\newblock {A note on Hilbert C*-modules associated with a foliation}.
\newblock {\em Publ. Res. Inst. Math. Sci.}, {\bf 20}, (1984).
\newblock 97-106.

\bibitem{Yamagami:94}
S.~Yamagami.
\newblock {Modular theory for bimodules}.
\newblock {\em J. Funct. Anal.}, {\bf 125}, (1994).
\newblock 327-357.

\bibitem{Yang:87}
Y.~Yang.
\newblock {Numerical ranges of operators on Hilbert C*-modules}.
\newblock {\em Bull. Korean Math. Soc.}, {\bf 24}, (1987).
\newblock no. 1, p. 52 (Abstract of thesis).

\bibitem{Yang:84}
Youngoh Yang.
\newblock {A note on the numerical range of an operator}.
\newblock {\em Bull. Korean Math. Soc.}, {\bf 21}, (1984).
\newblock 27-30.

\bibitem{Zekri:89}
R.~Zekri.
\newblock {A new description of Kasparov's theory of C*-algebra extensions}.
\newblock {\em J. Func. Anal.}, {\bf 84}, (1989).
\newblock 441-471.

\bibitem{Zeller--Meier:91/2}
G.~Zeller-Meier.
\newblock {Noyaux positifs a valeurs dans une C*-algebre}.
\newblock preprint, Marseille, 1991.

\bibitem{Zeller--Meier:91/1}
G.~Zeller-Meier.
\newblock {Some remarks about C*--Hilbert spaces and Hilbert C*-modules}.
\newblock preprint, Marseille, 1991.

\bibitem{Zeller--Meier:87}
G.~Zeller-Meier.
\newblock {Some so far apparently unnoticed remarks on Hilbert C*-modules}.
\newblock Oberwolfach, Tagungsbericht 43/1987, pp.16-17.

\bibitem{Zettl:82/2}
H.~H. Zettl.
\newblock {Ideals in Hilbert modules and invariants under strong Morita
  equivalence of C*-algebras}.
\newblock {\em Arch. Math.}, {\bf 39}, (1982).
\newblock 69-77.

\bibitem{Zettl:82/1}
H.~H. Zettl.
\newblock {Strong Morita equivalence of C*-algebras preserves nuclearity}.
\newblock {\em Arch. Math.}, {\bf 38}, (1982).
\newblock 448-452.

\bibitem{Zettl:83}
H.~H. Zettl.
\newblock {A characterization of ternary rings of operators}.
\newblock {\em Adv. Math.}, {\bf 48}, (1983).
\newblock 117-143.

\bibitem{Zhang:88}
S.~Zhang.
\newblock {\em {On the structure of multiplier algebras}}.
\newblock PhD thesis, Purdue University, West Lafayette, U.S.A., 1988.

\bibitem{Zhang:89}
S.~Zhang.
\newblock {Stable isomorphisms of hereditary C*-subalgebras and stable
  equivalence of open projections}.
\newblock {\em Proc. Amer. Math. Soc.}, {\bf 105}, (1989).
\newblock 677-682.

\bibitem{Zhang:90}
S.~Zhang.
\newblock {Diagonalizing projections in multiplier algebras and in matrices
  over a C*-algebra}.
\newblock {\em Pacific J. Math.}, {\bf 145}, (1990).
\newblock 181-200.

\bibitem{Zhang:91/1}
S.~Zhang.
\newblock {Ideals of generalized Calkin algebras}.
\newblock {\em Contemp. Math.}, {\bf 120}, (1991).
\newblock 193-198.

\bibitem{Zhang:91/4}
S.~Zhang.
\newblock {K-theory, K-skeleton factorizations and bi-variable index
  Index($x,p$), (I,II,III)}.
\newblock preprint, 1991.

\bibitem{Zhang:91/3}
S.~Zhang.
\newblock {$K_1$-groups, quasidiagonality, and interpolation by multiplier
  projections}.
\newblock {\em Trans. Amer. Math. Soc.}, {\bf 325}, (1991).
\newblock 793-818.

\bibitem{Zhang:91/5}
S.~Zhang.
\newblock {On the homotopy type of the unitary group and the Grasmann space of
  purely infinite simple C*-algebras}.
\newblock preprint, 1991.

\bibitem{Zhang:91/2}
S.~Zhang.
\newblock {Problems on C*-algebras of real rank zero and their multiplier
  algebras}.
\newblock {\em Contemp. Math.}, {\bf 120}, (1991).
\newblock 199-203.

\bibitem{Zhang:91/6}
S.~Zhang.
\newblock {Rectifiable diameters of the Grasmann spaces of von Neumann algebras
  and certain C*-algebras}.
\newblock preprint, 1991.

\bibitem{Zhang:92}
S.~Zhang.
\newblock {Certain C*-algebras with real rank zero and their corona and
  multiplier algebras, I}.
\newblock {\em Pacific J. Math.}, {\bf 155}, (1992).
\newblock 169-197.

\bibitem{Zhang:92/2}
S.~Zhang.
\newblock {Certain C*-algebras with real rank zero and their corona and
  multiplier algebras, II}.
\newblock {\em K-Theory}, {\bf 6}, (1992).
\newblock 1-27.

\bibitem{Zhang:92/3}
S.~Zhang.
\newblock {On the exponential rank and exponential length of C*-algebras}.
\newblock {\em J. Operator Theory}, {\bf 28}, (1992).
\newblock 337-355.

\bibitem{Zhang:92/4}
S.~Zhang.
\newblock {Torsion of K-theory; bi-variable index and certain invariants of the
  essential commutant of ${\rm M}_n({\bf C})$, (I,II)}.
\newblock preprint, 1992.

\bibitem{Zhang:93}
S.~Zhang.
\newblock {Exponential rank and exponential length of operators on Hilbert
  C*-modules}.
\newblock {\em Annals Math.}, {\bf 137}, (1993).
\newblock 129-144.

\bibitem{Zhang:93/2}
S.~Zhang.
\newblock {Factorizations of invertible operators and K-theory of C*-algebras}.
\newblock {\em Bull. Amer. Math. Soc.}, {\bf 28}, (1993).
\newblock 75-83.

\bibitem{Zhang:93/3}
S.~Zhang.
\newblock {K-theory and a bivariable Fredholm index}.
\newblock {\em Contemp. Math.}, {\bf 148}, (1993).
\newblock 155-190.

\bibitem{Zhang:94/1}
S.~Zhang.
\newblock {K-theory and homotopy of certain groups and infinite Grassmann
  spaces associated with C*-algebras}.
\newblock {\em Int. J. Math.}, {\bf 5}, (1994).
\newblock 425-445.

\end{thebibliography}
\end{document}